\documentclass{amsart}
\usepackage{amsmath}
\usepackage{amssymb}
\usepackage{color}
\usepackage{graphicx}

\newtheorem{theorem}{Theorem}[section]
\newtheorem{lemma}[theorem]{Lemma}

\theoremstyle{definition}
\newtheorem{definition}[theorem]{Definition}

\newtheorem{assumption}{Assumption}

\theoremstyle{remark}

\numberwithin{equation}{section}

\newcommand{\R}{\mathbb{R}}

\renewcommand{\S}{\mathbb{S}}
\newcommand{\bP}{\mathbb{P}}

\newcommand{\pa}{\partial}
\newcommand{\F}{\mathcal{F}}
\newcommand{\ve}{\varepsilon}
\DeclareMathOperator*{\dv}{div}
\newcommand{\E}{\mathbb{E}}
\renewcommand{\u}{{\bf u}}
\renewcommand{\d}{{\bf d}}
\allowdisplaybreaks
\newcommand{\h}{{\bf h}}
\newcommand{\bu}{\overline{\u}}
\newcommand{\bd}{\overline{\d}}
\newcommand{\bW}{\overline{W}}
\newcommand{\bB}{\overline{B}}
\newcommand{\f}{{\bf f}}
\DeclareMathOperator{\dist}{dist}
\newcommand{\x}{x}
\allowdisplaybreaks
\begin{document}

\title[Stochastic Ericksen-Leslie system]{Global weak solutions to the Stochastic Ericksen--Leslie equations in dimension two}

\author{Hengrong Du}
\address{Department of Mathematics, Purdue University, West Lafayette, Indiana 47906}
\email{du155@purdue.edu}

\author{Changyou Wang}
\address{Department of Mathematics, Purdue University, West Lafayette, Indiana 47906}
\email{wang2482@purdue.edu}
\thanks{Both authors are partially supported by NSF DMS 1764417.}

\subjclass[2000]{Primary ; Secondary }



\keywords{Stochastic Ericksen--Leslie equations, }

\begin{abstract}
We establish the global existence of weak martingale solutions to the simplified stochastic Ericksen--Leslie system modeling the nematic liquid crystal flow driven by Wiener-type noises on the two-dimensional bounded domains. The construction of  solutions is based on the convergence of Ginzburg--Landau approximations. To achieve such a  convergence, we first utilize the concentration-cancellation method for the Ericksen stress tensor fields based on a Pohozaev type argument, and second the Skorokhod compactness theorem, which is built upon a uniform energy estimate.    
\end{abstract}

\maketitle


\section{Introduction}
In this article, we consider the following simplified stochastic Ericksen--Leslie system on a two dimensional bounded domain $D$ with smooth boundary:
\begin{equation}
  \left\{
  \begin{array}{l}
    d \u+(\u\cdot \nabla \u+\nabla P-\mu\Delta \u)dt= -\lambda\nabla\cdot (\nabla \d\odot \nabla \d)dt+\xi_1S(\u)dW_1, \\
    \nabla\cdot  \u=0, \\
    d \d+\u\cdot \nabla \d dt=\gamma(\Delta \d+|\nabla \d|^2 \d )dt+ \xi_2(\d\times \h)\circ dW_2,
  \end{array}
  \right.
  \label{eqn:SEL}
\end{equation}
where $\u:D\times\R_+\times \Omega\to\R^2$, $\d:D\times\R_+\times \Omega\to \S^2$ represent the fluid velocity field and the molecular director field, respectively, $P:D\times \R_+\times \Omega\to \R$ stands for the hydro-static pressure. $(\nabla \d\odot \nabla \d)_{ij}=\langle\pa_i \d, \pa_j\d\rangle$ $(1\le i, j\le 2)$ represents the Ericksen stress tensor field. The multiplicative noise term $S(\u)dW_1$ in \eqref{eqn:SEL}$_1$ shall be understood in the It\^{o} sense with a cylindrical Wiener process $W_1$ on a separable Hilbert space $K_1$. For a given $\h:\R^2\to \R^3$, $(\d\times\h)\circ dW_2$ is understood in the Stratonovich sense with a standard real-valued Brownian motion $W_2$. 
$\mu, \lambda, \gamma, \xi_1, \xi_2$ are positive physical constants.

 We assume, further, $(\u, \d)$ satisfies the
following initial-boundary conditions:
\begin{equation}
  (\u, \d)|_{t=0}=(\u_0, \d_0),  \quad \text{ in }D.
  \label{eqn:InitialData}
\end{equation}
\begin{equation}
	\u|_{\pa D}=0, \quad \frac{\pa \d}{\pa {\bf n}}\bigg|_{\pa D}=0, \quad (\text{or } \d|_{\pa D}=\d_0). 
	\label{eqn:boundaryData}
\end{equation}
where ${\bf n}$ is the unit outward normal to $\pa D$. 
In this paper, we use the Ginzburg--Landau type approximation which relaxes the condition $|\d|=1$ in \eqref{eqn:SEL} by introducing a penalized term, more specifically, we have a family of solutions $(\u^\ve, \d^\ve)_{0<\ve<1}$ to  
\begin{equation}
  \left\{
  \begin{array}{l}
    d\u^\ve+(\u^\ve\cdot \nabla \u^\ve+\nabla P^\ve-\mu\Delta \u^\ve)dt\\
   =-\lambda\nabla\cdot (\nabla \d^\ve\odot \nabla \d^\ve)dt+\xi_1S(\u^\ve) dW_1, \\
    \nabla \cdot \u^\ve=0, \\
    d \d^\ve+\u^\ve\cdot \nabla \d^\ve dt=\gamma\Big(\Delta \d^\ve-\f_\ve(\d^\ve)\Big)dt+ \xi_2(\d^\ve\times \h) \circ dW_2, 
  \end{array}
  \right.
  \label{eqn:GZapprx}
\end{equation}
where $\f_\ve(\d^\ve)=\nabla_{\d}F_\ve(\d^\ve)=\frac{1}{\ve^2}(|\d^\ve|^2-1)\d^\ve$ with $F_\ve(\d)=\frac{1}{4\ve^2}(1-|\d|^2)^2$. 

In the deterministic case $(\xi_1=\xi_2=0)$, the global existence of the weak solutions to the Ginzburg--Landau type Ericksen--Leslie system \eqref{eqn:GZapprx} was first  investigated by Lin--Liu \cite{linliu1995nonparabolic} which is a simplified version of the full Ericksen--Leslie system \cite{ericksen1961conservation,Ericksen1962,Leslie1968,leslie1992continuum}. For the simplified Ericksen--Leslie system \eqref{eqn:SEL}, motivated by Struwe \cite{Struwe85} on harmonic map heat flows in dimension two,
the existence of a unique global weak solution with partial regularity was established Lin--Lin--Wang \cite{lin2010liquid} and Lin-Wang \cite{lin2010},
which was generalized  by Huang--Lin--Wang \cite{HuangLinWang2014EricksenLeslieTwoDimension} for the full Ericksen--Leslie system. 
See also Hong \cite{Hong2011} and Hong--Xin \cite{HongXin2011} for related works. 
We refer the readers to \cite{lin2014recent}
for a comprehensive survey for the recent developments. The question that whether or not one can obtain a weak solution of \eqref{eqn:SEL} via sending $\ve\to0$ in \eqref{eqn:GZapprx} remains open due to the difficulty with possible defect measures appearing in the Ericksen stress tensor field.  In a very recent paper \cite{kortum2019concentrationcancellation}, Kortum applied a concentration-cancellation method initiated by  Diperna--Majda \cite{DipernaMajda1988concentration} on 2-D incompressible Euler equation  to show that ${\rm div}(\nabla\d^\ve\odot\nabla\d^\ve)\rightharpoonup{\rm div}(\nabla\d\odot\nabla\d)$ in the torus $\mathbb{T}^2$. For general domains and full Ericksen--Leslie system, the weak compactness result was shown in \cite{DuHuangWang2020compactness} via the Hopf differential and the Pohozaev technique. We also want to point out that in 3-D, the Ginzburg--Landau approximation was implemented in \cite{linwang2016weakthreedim} to construct a global weak solutions to the simplified Ericksen--Leslie system \eqref{eqn:SEL} with the half-sphere assumption imposed on directors ($\d\in \mathbb{S}^2_+$).

On the other hand, there is a growing number of research studies devoted to the simplified stochastic Ericksen--Leslie system \eqref{eqn:GZapprx} with various types of random noises ($\xi_1^2+\xi_2^2> 0$). See for instance, \cite{brzezniak20202d,Brzezniakhausenblas2019NoteonSEL,brzezniak2019some,Brzeniak2019SELwithjumpnoise}. For the mathematically modeling, taking the stochastic terms into account reflects the influence of environmental noises,  measurement uncertainties or thermal fluctuations. Analogously, Bouard--Hocquet--Prohl obtained the Struwe-like global solution to \eqref{eqn:SEL} in  \cite{bouard2019existence} by a bootstrap argument together with Gy\"{o}ngy--Krylov $L^p$ estimates \cite{gyongy1996existence}. Very recently, 
Brze\'{z}niak, Deugou\'{e}, and Razafimandimby in \cite{brzezniak2020GL} proved the existence of short
time strong solutions to the simplified stochastic Ericksen--Leslie system. 
The main goal of this paper is to obtain a global weak solution to \eqref{eqn:SEL} by extending the compactness argument from \cite{DuHuangWang2020compactness} into the stochastic setting. 

For simplicity, we assume $\lambda=\xi_1=\gamma=\xi_2=1$. We introduce the notations of some function spaces:
\begin{align*}
    &{\bf H}=\text{closure of }C_0^\infty(D, \R^2)\cap\left\{ f|\nabla \cdot  f=0 \right\} \text{ in }L^2(D, \R^2),\\
    &{\bf J}=\text{closure of }C_0^\infty(D, \R^2)\cap\left\{ f|\nabla \cdot f=0 \right\}\text{ in }H_0^1(D, \R^2), \\
   & H^1(D, \S^2)=\left\{ f\in H^1(D, \R^3)| |f|=1 \text{ a.e. }x\in D \right\}. 
\end{align*}

 For a complete probability space $(\Omega, \F, \bP)$ with a filtration $\{\F_t\}_{t\ge 0}$, let $K_1$ be an infinite dimensional separable Hilbert space and  $W_1=\{W_1(t)\}_{t\ge 0}$ be a $K_1$-cylindrical Wiener process such that it is formally written as a series
 \begin{equation*}
 	W_1(t)=\sum_{i=1}^{\infty} B_i(t)e_i, \forall t\ge0,
 \end{equation*}
 where $\{B_i(t)\}_{i=1}^\infty$ is a family of i.i.d. standard Brownian motions and $\{e_i\}_{i=1}^\infty$ is an orthonormal base of $K_1$. The above series does not converge in $K_1$, but it does converge in $K_2$ if $K_2$ is a larger Hilbert space containing $K_1$ such that the inclusion map $J:K_1\to K_2$ is Hilbert-Schmidt. It is always possible to construct construct a space $K_2$ with this property. For example, we can define $K_2$ to be the closure of $K_1$ under the norm 
  \begin{equation*}
  \|x\|_{K_2}^2=\sum_{i=1}^{\infty}\frac{1}{i^2}\langle x, e_i \rangle_{K_1}^2.
  \end{equation*}
Then we can view $W_1$ as a $K_2$-valued Wiener process. Let $W_2=\{ W_2(t)\}_{t\ge 0}$ be a standard Brownian motion on $(\Omega, \F, \bP)$ adapted to $\{\F_t\}_{t\ge0}$. $S$ is a map from ${\bf H}$ to $\mathcal{L}_2(K_1, {\bf J})$, where $\mathcal{L}_2(K_1, {\bf J})$ denotes the  space of all Hilbert–Schmidt operators
from $K_1$ to ${\bf J}$, i.e., $\sum_{i=1}^{\infty}\|S(\cdot)(e_i)\|_{ {\bf J}}^2<\infty$, if $\{e_i\}_{i=1}^\infty$ is an orthonormal base of $K_1$. 


We now introduce the notion of a weak martingale solution to \eqref{eqn:SEL}. 
\begin{definition}\label{def:marsol}
A weak martingale solution to \eqref{eqn:SEL}, \eqref{eqn:InitialData}, \eqref{eqn:boundaryData} is a system consisting of a complete filtered probability space $(\Omega, \F, \bP)$ with a filtration $\{\F_t\}_{t\ge0}$, and $\F_t$ adapted stochastic processes $(\u(t), \d(t), W_1(t), W_2(t))_{t\ge0}$ such that for any $0<T<\infty$
  \begin{enumerate}
    \item $\{W_1(t)\}_{t\ge 0}$ (or $\{W_2(t)\}_{t\ge 0}$) is a $K_1$-cylindrical (resp. real-valued) Wiener process. 
    \item $(\u, \d):\Omega\times\R_+\to {\bf H}\times H^1(D, \S^2)$ is progressively measurable with respect to the filtration $\left\{ \mathcal{F}_t \right\}_{t\ge 0}$ such that for almost surely $\omega\in \Omega$, 
      \begin{equation*}
        \u\in L_t^\infty([0, T], {\bf H})\cap L_t^2([0,T], {\bf J}), \quad \d\in L^2([0,T], H^1(D, \S^2)).
      \end{equation*}
    \item We have
      \begin{equation}
        \begin{split}
          \E&\left[ \sup_{0\le t\le T}\int_{D\times\{t\}} |\u|^2+|\nabla \d|^2 +\int_{0}^{T}\int_D(|\nabla \u|^2+|\Delta \d+|\nabla \d|^2 \d|^2)dxds\right]\\
          &<\infty.
        \end{split}
        \label{eqn:defenergy}
      \end{equation}
      
         \item For almost surely $\omega\in \Omega$, for every $t\in [0,T]$, for any $\varphi\in C^\infty(D, \R^2)$, $\dv \varphi=0$, we have
      \begin{equation}
        \begin{split}
          &-\int_{D\times \left\{ t \right\}} \langle \u, \varphi \rangle dx-\int_{0}^{t}\int_{D}( \langle\u\otimes\u,\nabla \varphi\rangle +\langle  \u, \Delta \varphi \rangle)dxds\\
          &=-\int_{D}\left\langle \u_0, \varphi \right\rangle dx+\int_{0}^{t}\int_{D}(\langle \nabla \d\odot \nabla \d-\frac{1}{2}|\nabla \d|^2 \mathbb{I}_2, \nabla \varphi \rangle)dxds\\
          &+\int_{0}^{t}\int_D \langle \varphi, S(\u) dW_1(s) \rangle dx,
        \end{split}
        \label{eqn:uweak}
      \end{equation}
      and for any $\psi\in C^\infty(D, \R^3)$, 
      \begin{equation}
        \begin{split}
          &-\int_{D\times\{t\}}\langle \d, \psi \rangle dx-\int_{0}^{t}\int_{D}(\langle\u\otimes\d,\nabla \psi\rangle+\langle  \d, \Delta \psi \rangle)dxds\\
          &=-\int_{D}\langle  \d_0, \psi \rangle dx+\int_{0}^{t}\int_{D}\langle |\nabla \d|^2 \d, \psi \rangle dxds\\
          &+\int_{0}^{t}\int_{D}\langle\psi, (\d\times\h) \rangle dx \circ dW_2(s). 
        \end{split}
        \label{eqn:dweak}
      \end{equation}
        \end{enumerate}
\end{definition}
We introduce the following assumptions needed in our theorem. 
\begin{assumption}\label{assum}
Let $S:{\bf H}\to \mathcal{L}_2(K_1, {\bf J})$  be a global Lipschitz map. In particular, there exists $C>0$ such that $\|S(\u)\|^2_{\mathcal{L}_2(K_1, {\bf J})}\le
C(1+\|\u\|^2_{ {\bf H}})$ for all $\u\in {\bf H}$. $\h\in H^2(\R^2, \R^3)$. $(\u_0,\d_0)\in {\bf H}\times H^1(D;\mathbb{S}^2)$. Furthermore,  we assume $\{(\u_0^\ve, \d_0^\ve)\}_{0<\ve<1}\subset{\bf J}\times H^2(D;\mathbb{S}^2)$ and satisfies $(\u^\ve_0, \d_0^\ve)\to (\u_0, \d_0)$ in ${\bf H}\times H^1(D;\R^3)$.
\end{assumption}

Similar to Definition \ref{def:marsol}, a weak martingale solution $(\u^\ve(t), \d^\ve(t), W^\ve_1(t), W_2^\ve(t))$ adapted to a family of complete filtered probability spaces $(\Omega^\ve, \F^\ve, \bP^\ve, \{\mathcal{F}^\ve_t\}_{t\ge0})$ to \eqref{eqn:GZapprx}, \eqref{eqn:InitialData}, \eqref{eqn:boundaryData} can be defined. 
Under Assumption \ref{assum}, the existence of weak martingale solutions $(\u^\ve, \d^\ve, W_1^\ve, W_2^\ve)$ with respect to $ (\Omega^\ve, \mathcal{F}^\ve, \bP^\ve,\{\mathcal{F}^{\ve}_t\}_{t\ge0})$ was established in \cite[Theorem 3.2]{brzezniak2019some} via the Faedo--Galerkin approximation and compactness methods, together with the pathwise uniqueness in 2-D \cite[Theorem 3.4]{brzezniak2019some}. 
It has been proved in the recent work \cite[Theorem 3.17]{brzezniak2020strong} that \eqref{eqn:GZapprx} possesses a unique strong solution, that is, given $(\Omega, \F, \bP, \{\F_t\}_{t\ge 0}, W_1, W_2)$, there exists a unique pair of stochastic processes $(\u^\ve, \d^\ve)$ which solves \eqref{eqn:GZapprx} with respect to $(\Omega, \F, \bP, \{\F_t\}_{t\ge 0}, W_1, W_2)$ for initial data $(\u_0^\ve,\d_0^\ve)\in{\bf J}\times H^2(D;\R^3)$. 

Our main result states that we can obtain a global weak martingale solution to \eqref{eqn:SEL} via passing the limit  of soloutions $(\u^\ve,\d^\ve)$ to \eqref{eqn:GZapprx}:
\begin{theorem} Under Assumption \ref{assum}, 
  there exist a completed filtered probability space $(\Omega', \F',\bP')$ and a sequence of weak martingale solutions $(\bu^\ve, \bd^\ve, \bW_1^\ve, \bW_2^\ve)$ to \eqref{eqn:GZapprx}, \eqref{eqn:InitialData}, \eqref{eqn:boundaryData} on $(\Omega', \F', \bP')$ and  a weak martingale solution $(\u, \d, W_1', W_2')$ to \eqref{eqn:SEL}, \eqref{eqn:InitialData}, \eqref{eqn:boundaryData} such that after passing to a subsequence,
  \begin{equation*}
    \bu^\ve\rightharpoonup \u \text{ in }L^2(\Omega';L^2([0,T], H^1(D))), \quad \bd^\ve \rightharpoonup \d \text{ in }L^2(\Omega';L^2([0,T], H^1(D)))
  \end{equation*}
as $\ve\to0$.
  \label{thm:main}
\end{theorem}
The paper is organized as follows. In section 2 we establish some uniform energy estimates for the approximation solution $(\u^\ve, \d^\ve)$ by It\^{o}'s formula. The convergence of the approximated system, in particular,  the Ericksen stress tensor field and martingale terms will be discussed in section 3. In Appendix A, we provide the computation of It\^{o}'s formula for two functionals of $\d$. 
\section{Uniform estimates on approximated solutions}
  In this section, we will derive an uniform energy estimate for \eqref{eqn:GZapprx}, \eqref{eqn:InitialData}, \eqref{eqn:boundaryData} via It\^{o}'s calculus. 
  
  For simplicity, we denote $\|\cdot\|:=\|\cdot\|_{L^2(D)}$. 
   First,  applying It\^{o}'s formula to $\frac{1}{2}\|\u^\ve(t)\|^2$ yields 
\begin{align}  \label{eqn:Ito1}
    &\frac{1}{2}\|\u^\ve(t)\|^2-\frac{1}{2}\|\u_0^\ve\|^2+\int_{0}^{t}\int_D |\nabla\u^\ve|^2 dxds\\
    &=\int_{0}^{t}\int_D \langle \nabla \d^\ve\odot\nabla \d^\ve, \nabla \u^\ve \rangle dxds\nonumber\\
    & +\frac{1}{2}\int_{0}^{t}\left\|S(\u^\ve)\right\|_{\mathcal{L}_2(K_1, {\bf H})}^2 ds+\int_{0}^{t}\int_D\langle \u^\ve, S(\u^\ve)dW_1(s)\rangle dx ,  \nonumber
\end{align}
where we use the cancellation 
\begin{equation*}
  \int_{0}^{t}\int_D \langle \u^\ve\cdot \nabla \u^\ve, \u^\ve \rangle dxds=0. 
\end{equation*}
From the relation between Stratonovich and It\^{o}'s integral, we have that 
\begin{equation*}
  (\d\times \h)\circ dW_2=\frac{1}{2}((\d\times\h)\times\h) dt+(\d\times\h) dW_2. 
\end{equation*}
Therefore \eqref{eqn:GZapprx}$_3$ can be written as 
\begin{equation}
  d\d^\ve+\u^\ve\cdot \nabla \d^\ve dt=\left( \Delta\d^\ve-\f_\ve(\d^\ve) +\frac{1}{2}(\d^\ve\times\h)\times\h\right) dt+(\d^\ve\times\h) dW_2.
  \label{}
\end{equation}
Now we apply the It\^{o} formula to  $\Phi_\ve(\d^\ve):=\frac{1}{2}\|\nabla \d^\ve\|^2+\int_D F_\ve(\d^\ve)dx$ (see Appendix A) to get 
\begin{align}  \label{eqn:Ito2}
   	&\Phi_\ve(\d^\ve)(t)-\Phi_\ve(\d_0^\ve)\\
   	&=\int_{0}^{t}\int_D \langle \u^\ve\cdot \nabla\d^\ve, \Delta\d^\ve-\f_\ve(\d^\ve) \rangle dxds
	-\int_{0}^{t}\int_D |\Delta\d^\ve-\f_\ve(\d^\ve)|^2 dxds\nonumber\\
   	&+\frac{1}{2}\int_{0}^{t}\int_D (\langle \nabla\d^\ve, \nabla( (\d^\ve\times\h)\times\h) \rangle+|\nabla(\d^\ve\times\h)|^2) dxds\nonumber\\
   	&+\frac{1}{2}\int_{0}^{t}\int_D \langle -\Delta\d^\ve+\f_\ve(\d^\ve), \d^\ve\times\h \rangle dx dW_2(s).\nonumber
\end{align}
Using the fact that 
\begin{equation*}
\int_{0}^{t}\int_D \langle \u^\ve\cdot \nabla\d^\ve, \f_\ve(\d^\ve) \rangle dxds=\int_{0}^{t}\int_D \u^\ve\cdot \nabla F_\ve(\d^\ve)dxds=0,
\end{equation*}
and
\begin{align*}
    &-\int_{0}^{t}\int_D \langle \u^\ve \cdot \nabla \d^\ve, \Delta\d^\ve \rangle dxds\\
    &=\int_{0}^{t}\int_D \langle \nabla \d^\ve\odot \nabla\d^\ve, \nabla\u^\ve \rangle dxds+\int_{0}^{t}\int_D \u^\ve\cdot \nabla\left( \frac{|\nabla\d^\ve|^2}{2} \right) dxds\\
    &=\int_{0}^{t}\int_D  \langle \nabla\d^\ve\odot \nabla\d^\ve, \nabla\u^\ve \rangle dxds,
\end{align*}
we can add \eqref{eqn:Ito1} and \eqref{eqn:Ito2} together to obtain
\begin{align}  \label{eqn:EnergyId}
    &\frac{1}{2}\|\u^\ve(t)\|^2+\frac{1}{2}\|\nabla\d^\ve(t)\|^2\\
    &+\int_{D\times\{t\}}F_\ve(\d^\ve)dx+\int_{0}^{t}\int_D (|\nabla\u^\ve|^2 +|\Delta\d^\ve-\f_\ve(\d^\ve)|^2) dxds\nonumber\\
    &=\frac{1}{2}\|\u_0^\ve\|^2+\frac{1}{2}\|\nabla\d_0^\ve\|^2\nonumber\\
    &+\frac{1}{2}\int_{0}^{t}(\left\|S(\u^\ve)\right\|_{\mathcal{L}_2(K_1, {\bf H})}^2 +\|\nabla(\d^\ve\times\h)\|^2) ds\nonumber\\
    &+ \frac{1}{2}\int_{0}^{t}\int_D \langle \nabla\d^\ve,\nabla ((\d^\ve\times\h)\times\h ) \rangle dxds\nonumber\\
    &+\int_{0}^{t}\int_D\langle \u^\ve, S(\u^\ve)dW_1(s)  \rangle dx \nonumber\\
    &+ \int_{0}^{t}\int_D \langle \d^\ve\times\h, \Delta\d^\ve-\f_\ve(\d^\ve)\rangle dx dW_2(s). \nonumber
\end{align}
It has been shown in \cite[Theorem 5.1]{brzezniak2019some} that $\d^\ve$ satisfies the maximum principle, i.e., $|\d^\ve|\le 1$ for almost all $(\omega, t, x)\in \Omega\times[0,T]\times D$ provided $|\d_0^\ve|\le 1$.
 Hence we have that 
\begin{align*}
  &\int_{0}^{t}\left\|S(\u^\ve)\right\|_{\mathcal{L}_2(K_1, {\bf H})}^2 ds\le \int_{0}^{t}\left\|S(\u^\ve)\right\|_{\mathcal{L}_2(K_1, {\bf J})}^2ds \le C \int_{0}^{t}\int_D (1+|\u^\ve|^2) dxds, \\
  &\int_{0}^{t}\int_D|\nabla(\d^\ve\times\h)|^2 dxds\le C \int_{0}^{t}\int_D (|\nabla\d^\ve|^2+|\nabla\h|^2) dxds, \\
  &\int_0^t\int_D \langle \nabla \d^\ve, \nabla ((\d^\ve\times\h)\times\h)\rangle dxds\le C\int_0^t\int_D (|\nabla \d^\ve|^2+|\nabla \h|^2) dxds. 
\end{align*}
Combine all these estimates above, we arrive at 
\begin{align} \label{eqn:uniformineq}
  &\frac{1}{2}\|\u^\ve(t)\|^2+\frac{1}{2}\|\nabla\d^\ve(t)\|^2\\
  &+\int_{D\times\{t\}}F_\ve(\d^\ve)dx+\int_{0}^{t}\int_D (|\nabla\u^\ve|^2+|\Delta\d^\ve-\f_\ve(\d^\ve)|^2) dxds\nonumber\\
  &\le \frac{1}{2}\|\u_0^\ve\|^2+\frac{1}{2}\|\nabla\d_0^\ve\|^2\nonumber\\
  &+C\int_{0}^{t}\int_D (|\u^\ve|^2+|\nabla\d^\ve|^2+|\nabla\h|^2)dxds\nonumber\\
  &+\int_{0}^{t}\int_D\langle \u^\ve, S(\u^\ve) dW_1(s) \rangle dx +\int_{0}^{t}\int_D \langle \d^\ve\times\h, \Delta\d^\ve-\f_\ve(\d^\ve)\rangle dx dW_2(s) . \nonumber 
\end{align}
We can derive from taking the expectation of \eqref{eqn:uniformineq} that 
\begin{align}\label{eqn:Energyineq}
      &\E\sup_{0\le t\le T}\left[ \|\u^\ve(t)\|^2+\|\nabla\d^\ve(t)\|^2+\int_{D\times\{t\}}F_\ve(\d^\ve) dx \right]\\
      &+\E\int_{0}^{T}\int_D (|\nabla\u^\ve|^2+|\Delta\d^\ve-\f_\ve(\d^\ve)|^2) dxds\nonumber\\
      &\le C\E\int_{0}^{T}\int_D(|\u^\ve|^2+|\nabla\d^\ve|^2+|\nabla\h|^2 )dxds\nonumber\\
      & +C\E \sup_{0\le t\le T} \left|\int_{0}^{t}\int_D \langle \u^\ve(s), S(\u^\ve(s)) dW_1(s) \rangle dx \right|\nonumber\\
      & + C\E \sup_{0\le t\le T}\left|\int_{0}^{t}\int_D \langle \d^\ve\times\h, \Delta\d^\ve-\f_\ve(\d^\ve) \rangle dx dW_2(s)\right|\nonumber\\
      &+ C(1+\left\|\u_0\right\|^2+\|\nabla\d_0\|^2).\nonumber
\end{align}
Now we use the Burkholder--Davis--Gundy inequality, Cauchy--Schwarz  inequality and H\"{o}lder inequality to show that 
\begin{align}  \label{eqn:BDG1}
    &\E \sup_{0\le t\le T}\left|\int_{0}^{t}\int_D\langle \d^\ve\times\h, \Delta\d^\ve-\f_\ve(\d^\ve) \rangle dx dW_2(s)\right|\\
    &\le C\E\left[ \int_{0}^{T}\left|\int_D \langle \d^\ve\times\h, \Delta\d^\ve-\f_\ve(\d^\ve) \rangle dx\right|^2 ds \right]^{\frac{1}{2}}\\
    &\le C \E\left[ \int_{0}^{T}\|\langle \d^\ve\times\h, \Delta\d^\ve-\f_\ve(\d^\ve) \rangle\|^2 ds \right]^{\frac{1}{2}}\nonumber\\
    & \le C \E\left[ \sup_{0\le t\le T}\|\d^\ve\times\h\|_{L^\infty(D)}\left( \int_{0}^{T}\int_D |\Delta\d^\ve-\f_\ve(\d^\ve)|^2 dxds \right)^{\frac{1}{2}} \right]\nonumber\\
    & \le C \E\sup_{0\le t\le T}\|\d^\ve\times\h\|_{L^\infty(D)}^2+\frac{1}{4}\E\int_{0}^{T}\int_D |\Delta\d^\ve-\f_\ve(\d^\ve)|^2 dxds\nonumber\\
    &\le C(\left\|\h\right\|_{L^\infty}, T, D)+ \frac{1}{4} \E\int_{0}^{T}\int_D |\Delta\d^\ve-\f_\ve(\d^\ve)|^2 dxds.\nonumber
\end{align}
Similarly, we can show 
\begin{align}  \label{eqn:BDG2}
    &\E \sup_{0\le t\le T}\left|\int_{0}^{t}\int_D\langle \u^\ve(s), S(\u^\ve(s)) \rangle dx dW_1(s)\right|\\
    & \le \frac{1}{4}\E\sup_{0\le t\le T}\|\u^\ve(t)\|^2+C\E\int_{0}^{T}\int_D |\u^\ve|^2 dxds. \nonumber
\end{align}
Now we can substitute \eqref{eqn:BDG1} and \eqref{eqn:BDG2} into \eqref{eqn:Energyineq} to get 
\begin{align*}
    &\E\sup_{0\le t\le T}\left[ \|\u^\ve\|^2+\|\nabla\d^\ve\|^2+\int_{D\times\{t\}}F_\ve(\d^\ve)dx\right]\\
    &+\E\int_{0}^{T}\int_D (|\nabla\u^\ve|^2+|\Delta\d^\ve-\f_\ve(\d^\ve)|^2) dxds\\
    &\le C \E\int_{0}^{T}\int_D (|\u^\ve|^2+|\nabla\d^\ve|^2+|\nabla\h|^2)dxds+C(\|(\u_0, \nabla\d_0)\|, \|\h\|_{L^\infty}, T, D). 
\end{align*}
It follows from  Gronwall's lemma that 
\begin{align}  \label{eqn:unifEnergy}
    &\E\sup_{0\le t\le T}\left[ \|\u^\ve(t)\|^2+\|\nabla\d^\ve(t)\|^2+ \int_{D\times\{t\}}F_\ve(\d^\ve)dx \right]\\
    &+\E\int_{0}^{T}\int_D (|\nabla\u^\ve|^2+|\Delta\d^\ve-\f_\ve(\d^\ve)|^2) dxds\nonumber\\
    &\le C(\|(\u_0, \nabla\d_0)\|, \|\h\|_{L^\infty}, \|\nabla\h\|, T, D). \nonumber
\end{align}
Furthermore, if we raise both sides of \eqref{eqn:uniformineq} to the power $p$ $(p>1)$ and take the expectation, we arrive at 
\begin{align}\label{eqn:Lp1}
  &\E\sup_{0\le t\le T}\left[ \|\u^\ve(t)\|^2+\|\nabla\d^\ve(t)\|^2+\int_{D\times\{t\}}F_\ve(\d^\ve)dx  \right]^p\\
  & +\E\left[ \int_0^t\int_D (|\nabla\u^\ve|^2+|\Delta\d^\ve-\f_\ve(\d^\ve)|^2) dxds \right]^p\nonumber\\
  &\le C(\|\u_0\|, \|\nabla\d_0\|, p)+CT \E \big(\int_{0}^{T} \left[\|\u^\ve(t)\|^2+\|\nabla\d^\ve(t)\|^2+\|\nabla\h\|^2\right]dt\big)^p\nonumber\\
  & +C \E\sup_{0\le t\le T}\left|\int_{0}^{t}\int_D\langle \u^\ve(s), S(\u^\ve(s)) \rangle dx dW_1(s) \right|^p\nonumber\\
  &+C\E\sup_{0\le t\le T}\left|\int_{0}^{t}\int_D \langle \d^\ve\times\h, \Delta\d^\ve-\f_\ve(\d^\ve) \rangle dx dW_2(s)\right|^p.\nonumber
\end{align}
Now we apply the Burkholder--Davis--Gundy, Cauchy--Schwarz,  and H\"{o}lder inequalities to the last two terms in the righg hand side  to get
\begin{align}\label{eqn:Lp2}
  &\E\sup_{0\le t\le T}\left|\int_{0}^{t}\int_D \langle \u^\ve(s), S(\u^\ve(s))dW_1(s) \rangle dx \right|^p\\
  &\le C\E\left[ \int_{0}^{T} \|\u^\ve(s)\|^2\|S(\u^\ve(s))\|^2 ds \right]^{\frac{p}{2}}\nonumber\\
  &\le C\E\left[ \sup_{0\le t\le T}\|\u^\ve(t)\|^{p}\left( \int_{0}^{T}(1+\|\u^\ve(s)\|^2)ds \right)^{\frac{p}{2}} \right]\nonumber\\
  &\le \frac{1}{4}\E\sup_{0\le t\le T}\|\u^\ve(t)\|^{2p}+C\E\int_{0}^{T}(1+\|\u^\ve(s)\|^2)^p ds.\nonumber
\end{align}
A similar argument yields 
\begin{align}\label{eqn:Lp3}
  &\E\sup_{0\le t\le T}\left|\int_{0}^{t}\int_D \langle \d^\ve\times\h, \Delta\d^\ve-\f_\ve(\d^\ve) \rangle dx dW_2(s)\right|^p\\
  &\le \frac{1}{4}\E \left[ \int_{0}^{T}\int_D |\Delta\d^\ve-\f_\ve(\d^\ve)|^2 dxds\right]^p+ C\E\sup_{0\le t\le T}\|\d^\ve\times\h\|^{2p}.\nonumber 
\end{align}
Combine \eqref{eqn:Lp1}, \eqref{eqn:Lp2} and \eqref{eqn:Lp3}, by Gronwall's inequality we obtain that for $p\ge 1$, it holds
\begin{align}  \label{eqn:unifpnorm}
    &\E\sup_{0\le t\le T}\left[ \|\u^\ve(t)\|^2+\|\nabla\d^\ve(t)\|^2+ \int_{D\times\{t\}}F_\ve(\d^\ve)dx \right]^p\\
    &+\E\left[\int_{0}^{T}\int_D (|\nabla\u^\ve|^2+|\Delta\d^\ve-\f_\ve(\d^\ve)|^2) dxds\right]^p\nonumber\\
    &\le C(\|(\u_0, \nabla\d_0)\|, \|\h\|_{L^\infty}, \|\nabla\h\|, T, D, p). \nonumber
\end{align}

Similar to the Aubin-Lions lemma in the deterministic case, we need some fractional Sobolev estimates in $t$ variable as in \cite{flandoli1995martsol} for stochastic Navier-Stokes equations. 
Write
\begin{align*}
  \u^\ve(t)&=\u_0^\ve+\int_{0}^{t}{\bf P}\Delta \u^\ve(s)ds-\int_{0}^{t}{\bf P}\nabla\cdot (\u^\ve\otimes\u^\ve)(s)ds\\
  &-\int_{0}^{t}{\bf P}\nabla\cdot (\nabla\d^\ve\odot \nabla \d^\ve)(s)ds+\int_{0}^{t}
  S(\u^\ve(s))dW_1(s)\\
  &:=\u_0^\ve+\sum_{i=1}^{4}I_i^\ve(t),
\end{align*}
where ${\bf P}$ is the Leray projection operator. 
We have that 
\begin{align*}
  &\E\left[ \|I_1^\ve\|^2_{W^{1,2}([0,T]; H^{-1}(D))}+\|I_2^\ve\|^2_{W^{1,2}([0,T]; H^{-1}(D))} \right]\le C,\\
  &\E\left[ \|I_3^\ve\|^2_{W^{1,2}([0,T]; W^{-2, \tilde{p}}(D))} \right]\le C,\text{ for some } \tilde{p}>2.
\end{align*}
Applying \cite[Lemma 2.1]{flandoli1995martsol} to $I_4^\ve$ we conclude that for any $\alpha\in(0,\frac{1}{2})$ and $p\in[2,\infty)$, it holds
\begin{align*}
  \E[\|I_4^\ve\|_{W^{\alpha, p}([0,T]; L^2(D))}^p]&=\E\left\|\int_{0}^{t}S(\u^\ve(s))dW_1(s)\right\|_{W^{\alpha, p}([0, T]; L^2(D))}^p\\
  &\le C\E\int_{0}^{T}\|S(\u^\ve(t))\|_{\mathcal{L}_2(K_1, {\bf H})}^pdt\\
  &\le C \E \int_{0}^{T}(1+\|\u^\ve(t)\|_{L^2(D)}^p)dt\le C.  
\end{align*}
Now we define
\begin{align*}
  X&:=L^\infty([0,T];L^2(D))\cap L^2([0,T];H^1(D))\\
  &\cap \left( W^{1, 2}([0,T]; H^{-1}(D))+W^{1, 2}([0,T]; W^{-2, \tilde{p}}(D))+W^{\alpha, p}([0,T]; L^2(D)) \right).
\end{align*}
Let $\{\mathcal{L}(\u^\ve)\}_{0<\ve<1}$ be a family of probability measures define on $X$ as following:
\begin{equation*}
  \mathcal{L}(\u^\ve)(B):=\bP(\u^\ve\in B)
\end{equation*}
for any Borel set $B\subset X$. 
For a fix $R>0$, we can derive from Chebyshev's inequality that
\begin{align*}
	&\bP(\|\u^\ve\|_{X}>R)\\
	&\le \bP\left( \|\u^\ve\|_{L^\infty([0,T]; L^2(D))}>\frac{R}{3} \right)+\bP\left( \|\u^\ve\|_{L^2([0,T];H^1(D))}>\frac{R}{3} \right)\\
	&+\bP\left( \|\u^\ve\|_{W^{1,2}([0,T];H^{-1}(D))+W^{1, 2}([0,T];W^{-2,p}(D))+W^{\alpha, p}([0,T];L^2(D))}>\frac{R}{3} \right)\\
	&\le \frac{C}{R}. 
\end{align*}
By a fractional version of Aubin-Lions lemma and the Sobolev interpolation inequality, $X$ is compactly embedded in $L^p([0,T]; L^{p}(D))\cap C([0,T]; W^{-2,\tilde{p}}(D))$ for $1<p<4$ (c.f. \cite{flandoli1995martsol,simon1990sobolev}).  Therefore $\{\mathcal{L}(\u^\ve)\}_{0<\ve<1}$ is tight in $L^p([0,T]; L^{p}(D))\cap C([0,T]; W^{-2,\tilde{p}}(D))$ for $1<p<4$.
Similarly, we have
\begin{align*}
  \d^\ve(t)&=\d_0^\ve-\int_{0}^{t}\nabla \cdot (\u^\ve\otimes \d^\ve)(s)ds+\int_{0}^{t}(\Delta\d^\ve-\f_\ve(\d^\ve))(s)ds\\
  &+\frac{1}{2}\int_{0}^{t}((\d^\ve\times\h)\times\h)(s) ds+\int_{0}^{t}(\d^\ve\times\h)(s)dW_2(s)\\
  &:=\d_0^\ve+\sum_{i=1}^{4}J_i^\ve(t).
\end{align*}
Then we have 
\begin{align*}
  &\E\left[ \|J^\ve_1\|_{W^{1,\frac{4}{3}}([0,T]; L^{\frac{4}{3}}(D))}^{\frac{4}{3}} \right]<C, \\
  &\E\left[ \|J^\ve_2\|_{W^{1,2}([0,T]; L^2(D))}^2+\|J_3^\ve\|_{W^{1,2}([0,T]; L^\infty(D))}^2 \right]\le C, 
\end{align*}
and by an argument similar to that of $I_4^\ve$ we can show that for any $\alpha\in(0,\frac{1}{2})$ and $p\in[2,\infty)$, it holds
\begin{align*}
  \E\left[ \left\|J_4^\ve\right\|_{W^{\alpha, p}([0,T]; L^2(D))}^p \right]&=\E\left\|\int_{0}^{t}\d^\ve\times\h dW_2(s)\right\|_{W^{\alpha, p}([0,T]; L^2(D))}^p\\
  &\le C\E\int_{0}^{T}\left\|\d^\ve\times\h(t)\right\|_{L^2(D)}^p dt\\
  &\le C\E\int_{0}^{T}\left\|\h\right\|_{L^\infty}^p \|\d^\ve(t)\|_{L^2(D)}^{p}dt\le C. 
\end{align*}
Hence, the laws $\{\mathcal{L}(\d^\ve)\}_{0<\ve<1}$ are bounded in probability in 
\begin{align*}
 &Y:= L^\infty([0,T]; H^1(D))
  \\
  &\cap \left( W^{1, \frac{4}{3}}([0,T]; L^{\frac{4}{3}}(D))+W^{1,2}([0,T]; L^2(D))+W^{\alpha, p}([0,T]; L^2(D)) \right).
\end{align*}
\begin{sloppypar}
Since $Y$ is compactly embedded into $L^q([0,T];L^q(D))\cap C([0,T]; L^{\frac{4}{3}}(D))$, $p>1$,  $\{\mathcal{L}(\d^\ve)\}_{0<\ve<1}$ is tight in $L^q([0,T];L^q(D))\cap C([0,T]; L^{\frac{4}{3}}(D))$, $p>1$. 
\end{sloppypar}
\section{Convergence of Ginzburg-Landau approximation}
The main purpose of this section is mainly devoted to show the convergence of Ericksen stress tensor and the martingale terms. 
From the uniform energy estimates in the previous section, we know that 
$(\mathcal{L}(\u^\ve), \mathcal{L}(\d^\ve))$ is tight in $L^p([0,T]; L^{p}(D))\cap C([0,T]; W^{-2,\tilde{p}}(D))\times L^q([0,T];L^q(D))\cap C([0,T]; L^{\frac{4}{3}}(D))$ for $1<p<4, 1<q<\infty$. %
Now we apply the Prohorov's theorem, there exists a probability measure $\mu$ on $L^p([0,T]; L^{p}(D))\cap C([0,T]; W^{-2,\tilde{p}}(D))\times L^q([0,T];L^q(D))\cap C([0,T]; L^{\frac{4}{3}}(D))\times C([0,T]; K_2)\times C([0,T];\R)$, $1<p<4, 1<q<\infty$ such that after passing to a subsequence,  
\begin{equation*}
  \mathcal{L}(\u^\ve, \d^\ve, W_1, W_2)\rightharpoonup \mu. 
\end{equation*}
Then by Skorokhod's embedding theorem, there exists a complete probability space $(\Omega', \mathcal{F}',\bP')$ and a sequence of random variables $(\bu^\ve, \bd^\ve, \bW_1^\ve, \bW_2^\ve)$ on $(\Omega', \mathcal{F}', \bP')$ such that 
\begin{equation}
  \mathcal{L}(\bu^\ve, \bd^\ve, \bW_1^\ve, \bW_2^\ve)=\mathcal{L}(\u^\ve, \d^\ve, W_1, W_2),\label{eqn:law}
\end{equation}
and $(\u, \d, W_1', W_2')$ defined on $(\Omega', \mathcal{F}', \bP')$ such that
\begin{equation}\left\{
  \begin{array}{l}
    \mathcal{L}(\u, \d, W_1', W_2')=\mu, \\
    \bu^\ve\to \u \text{ in }L^p([0,T]; L^{p}(D))\cap C([0,T]; W^{-2,\tilde{p}}(D)), 1<p<4, \bP'\text{-a.s.},\\
    \bu^\ve\rightharpoonup \u \text{ in }L^2(\Omega'\times[0,T]; {\bf J}), \\
    \bd^\ve \to \d \text{ in }L^q([0,T];L^q(D))\cap C([0,T]; L^{\frac{4}{3}}(D)),1<q<\infty, \bP'\text{-a.s.},\\
    \bd^\ve\rightharpoonup \d \text{ in }L^2(\Omega'\times[0,T], H^1), \\
    W_1^\ve\to W_1' \text{ in }C([0,T]; K_2), \quad \bP'\text{-a.s.}, \\
    W_2^\ve\to W_2' \text{ in }C([0,T];\R), \quad \bP'\text{-a.s.}. 
  \end{array}
  \right.
  \label{eqn:conv}
\end{equation}
And for $\bP'$-a.s.,  $\u\in L^\infty([0,T]; {\bf H})\cap L^2([0,T]; {\bf J})$, $\d\in L^\infty([0,T]; H^1(D))$. 

For martingale solutions, for each $0<\ve<1$, we define $M_{\u^\ve}(t), M_{\d^\ve}(t)$ as
\begin{align*}
  M_{\u^\ve}(t)&=\u^\ve(t)-\u_0^\ve+\int_{0}^{t}[{\bf P}\nabla \cdot (\u^\ve\otimes \u^\ve)-{\bf P}\Delta \u^\ve+{\bf P}\nabla\cdot (\nabla\d^\ve\odot\nabla\d^\ve)](s)ds, \\
  M_{\d^\ve}(t)&=\d^\ve(t)-\d_0^\ve+\int_{0}^{t}[\nabla \cdot(\u^\ve\otimes\d^\ve)-\Delta \d^\ve+\f_\ve(\d^\ve)-\frac{1}{2}(\d^\ve\times\h)\times\h](s)ds,
 \label{}
\end{align*}
for any $t\in(0,T]$. Also define $M_{\bu^\ve}, M_{\bd^\ve}$ by replacing $\u^\ve, \d^\ve$ in $M_{\u^\ve}, M_{\d^\ve}$ by $\bu^\ve, \bd^\ve$. 

Next we show that for $\bP'$-a.s.,
\begin{align}
	M_{\bu^\ve}(t)&=\int_{0}^{t}S(\bu^\ve)d\bW_1^\ve(s), \label{eqn:Mbu}\\
	M_{\bd^\ve}(t)&=\int_{0}^{t}(\bd^\ve\times\h)d\bW_2^\ve(s)\label{eqn:Mbd}
\end{align}
for every $\ve>0$ and every $t\in[0,T]$.  For any ${\bf z}\in L^2(0,T; H^{-1})$ we set
\begin{equation*}
	\varphi({\bf z})=\frac{\int_{0}^{T}\|{\bf z}(s)\|_{H^{-1}}^2 ds}{1+\int_{0}^{T}\|{\bf z}(s)\|_{H^{-1}}^2 ds}.
\end{equation*}
By a argument similar to  that in \cite{bensoussan1995stochastic,brzezniak2019some} we can show that
\begin{equation*}
	\E' \varphi\left( M_{\bu^\ve}(\cdot )-\int_{0}^{\cdot }S(\bu^\ve(s))d\bW_1^\ve(s) \right)=\E\varphi\left( M_{\u^\ve}(\cdot )-\int_{0}^{\cdot }S(\u^\ve(s))dW_1(s)\right)=0.
\end{equation*}
This implies that for $\bP'$-a.s. \eqref{eqn:Mbu} holds for all $t\in(0,T]$. Similarly, we can show \eqref{eqn:Mbd} is also true. 

 Let $M_{\u}(t)$ and $M_{\d}(t)$ be defined by 
  \begin{align*}
    M_{\u}(t)&=\u(t)-\u_0+\int_{0}^{t}[{\bf P}\nabla\cdot (\u\otimes\u)-{\bf P}\Delta\u+{\bf P}\nabla \cdot (\nabla\d\odot\nabla\d)](s)ds, \\
    M_{\d}(t)&=\d(t)-\d_0+\int_{0}^{t}[\nabla\cdot (\u\otimes\d)-\Delta\d-|\nabla\d|^2\d-\frac{1}{2}(\d\times\h)\times\h](s)ds.
  \end{align*}
 With \eqref{eqn:conv}, we have the almost surely convergence of every term in $M_{\bu^\ve}$ except the Ericksen stress tensor $(\nabla\bd^\ve\odot\nabla\bd^\ve)$. Now we claim that for $\bP'$-a.s. 
\begin{align}\label{eqn:convofEricksen}
  &\lim_{\ve\to0}\int_{0}^{T}\int_D \langle \nabla\bd^\ve\otimes\nabla\bd^\ve-\frac{1}{2}|\nabla \bd^\ve|^2\mathbb{I}_2, \nabla \varphi \rangle dxds\\
  &=\int_{0}^{T}\int_D \langle \nabla \d\otimes\nabla\d-\frac{1}{2}|\nabla \d|^2\mathbb{I}_2, \nabla\varphi \rangle dxds. \nonumber
\end{align}
For any $0<\Lambda_1, \Lambda_2<\infty$, define the set ${\bf X}(\Lambda_1, \Lambda_2)$ consisting of solutions $\bd^\ve$ to
\begin{equation}
  \Delta \bd^\ve-\f_\ve(\bd^\ve)=\tau^\ve \text{ in }D
  \label{}
\end{equation}
such that the following properties hold:
\begin{enumerate}
  \item $|\bd^\ve|\le 1$ for a.e. $x\in D$. 
  \item $$\sup_{0<\ve\le 1} \mathcal{E}_\ve(\bd^\ve)=\int_D \left(\frac{1}{2}|\nabla\bd^\ve|^2+F_\ve(\bd^\ve)\right) dx\le \Lambda_1. $$
  \item \begin{equation*}
      \sup_{0<\ve\le 1}\|\tau^\ve\|_{L^2(D)}\le \Lambda_2.
    \end{equation*}
\end{enumerate}
The following small energy regularity lemma \cite{kortum2019concentrationcancellation,linwang2016weakthreedim} plays a key role in our analysis. 
\begin{lemma}
  Suppose $\{\bd^\ve\}_{0<\ve\le 1} \subset{\bf X}(\Lambda_1, \Lambda_2)$ and $\tau^\ve\rightharpoonup \tau$ in $L^2(D)$. Then there exists a $\delta_0>0$ such that if for $x_0\in D$ and $0<r_0<\dist(x_0, \pa\Omega)$,
  \begin{equation}
  \sup_{0<\ve\le 1}\int_{B_{r_0}(x_0)}\left(\frac{1}{2}|\nabla\bd^\ve|^2+F_\ve(\bd^\ve)\right)dx\le \delta_0^2,
\end{equation}
then there exists an approximated harmonic map $\d\in H^1(B_{\frac{r_0}{4}}(x_0), \mathbb{S}^2)$ with tensor field $\tau$, i.e., 
\begin{equation}
  \Delta\d+|\nabla\d|^2\d=\tau,
  \label{}
\end{equation}
such that 
\begin{equation}
  \bd^\ve\to \d \text{ in }H^1(B_{\frac{r_0}{4}}(x_0)) 
  \label{}
\end{equation}
as $\ve\to 0$. 
  \label{lemma:smallenergy}
\end{lemma}
This leads to the following $H^1$ precompactness result.
\begin{lemma}
  Under the same assumption as Lemma \ref{lemma:smallenergy},
  \begin{equation*}
    \bd^\ve\to \d \text{ in }H^1_{\text{loc}}(D\setminus \Sigma),
  \end{equation*}
  where 
  \begin{equation*}
    \Sigma:=\bigcap_{r>0}\left\{ x\in D:\liminf_{\ve\to0}\int_{B_r(x)}\left(\frac{1}{2}|\nabla\bd^\ve|^2+F_\ve(\bd^\ve)\right)dx>\delta_0^2 \right\}.
  \end{equation*}
  Moreover, $\Sigma$ is a finite set.
  \label{lemma:H1precom}
\end{lemma}
From \eqref{eqn:unifEnergy} and $\eqref{eqn:law}$, we have
\begin{align}\label{eqn:globeng}
  &\E' \sup_{0\le t\le T}\left[ \|\bu^\ve(t)\|^2+\|\nabla\bd^\ve(t)\|^2+\int_{D\times\{t\}} F_\ve(\bd^\ve) dx \right] \\
  & +\E'\left[\int_{0}^{T}(\|\nabla\bd^\ve\|^2+\|\Delta\bd^\ve-\f_\ve(\bd^\ve)\|^2)dt \right]\nonumber\\
   &=\E \sup_{0\le t\le T}\left[ \|\u^\ve(t)\|^2+\|\nabla\d^\ve(t)\|^2+\int_{D\times\{t\}} F_\ve(\d^\ve) dx\right] \nonumber\\
  & +\E\left[\int_{0}^{T}(\|\nabla\d^\ve\|^2+\|\Delta\d^\ve-\f_\ve(\d^\ve)\|^2)dt\right]\nonumber\\
  &\le C. \nonumber
\end{align}
Hence, there exists $\mathcal{N}\subset \Omega'$ such that $\bP'(\mathcal{N})=0$, and it holds for $\omega\in \Omega'\setminus \mathcal{N}$ that  
\begin{equation}
  \liminf_{\ve\to0} \int_{0}^{T}\int_D (|\nabla\bd^\ve|^2+|\Delta\bd^\ve-\f_\ve(\bd^\ve)|^2)dxdt =C_1(\omega)<\infty,
\end{equation}
and
\begin{equation}
  \liminf_{\ve\to 0}\sup_{0\le t\le T} \int_{D\times\{t\}}(|\bu^\ve|^2+|\nabla\bd^\ve|^2+F_\ve(\bd^\ve))dx=C_2(\omega)<\infty.
  \label{}
\end{equation}
Now fix $\omega\in \Omega'\setminus \mathcal{N}$, by Fatou's lemma, we have
\begin{align*}
  &\int_{0}^{T}\liminf_{\ve\to 0}\int_D (|\nabla\bd^\ve|^2+|\Delta\bd^\ve-\f_\ve(\bd^\ve)|^2)dxds\\
  &\le \liminf_{\ve\to 0}\int_{0}^{T}\int_D (|\nabla\bd^\ve|^2+|\Delta\bd^\ve-\f_\ve(\bd^\ve)|^2)dxds<\infty.
  \label{}
\end{align*}
Hence there exists $A\subset[0,T]$ with full Lebesgue such that for any $t\in A$, 
\begin{equation*}
  \liminf_{\ve\to 0}\int_{D\times\{t\}} (|\nabla\bd^\ve|^2+|\Delta\bd^\ve-\f_\ve(\bd^\ve)|^2)dx<\infty.
\end{equation*}
For $t\in A$, we set
\begin{equation*}
  \Sigma_t:=\bigcap_{r>0}\left\{ x \in D:\liminf_{\ve\to 0}\int_{B_r(x)\times\{t\}}(\frac{1}{2}|\nabla\bd^\ve|^2+F_\ve(\bd^\ve))dx>\delta_0^2 \right\}.
\end{equation*}
By Lemma \ref{lemma:H1precom}, it holds that $\#(\Sigma_t)\le C_3(\omega)<\infty$ and
\begin{equation*}
  \bd^\ve(t)\to \d(t) \text{ in }H^1_{\text{loc}}(D\setminus \Sigma_t). 
\end{equation*}
Hence we get \eqref{eqn:convofEricksen} holds
for $\varphi$ with $\text{supp }\varphi\subset D\setminus \Sigma_t$. Now we consider the case $\Sigma_t\cap\text{supp }\varphi\neq \emptyset$. Since $\Sigma_t$ is finite, we may assume $(0,0)\in\text{supp }\varphi$. Write
\begin{equation}
  \nabla\bd^\ve\odot \nabla\bd^\ve-\frac{1}{2}|\nabla\bd^\ve|^2 \mathbb{I}_2=\frac{1}{2}\left( 
  \begin{matrix}
    |\pa_{x_1} \bd^\ve|^2-|\pa_{x_2}\bd^\ve|^2 & 2\langle \pa_{x_1} \bd^\ve, \pa_{x_2} \bd^\ve\rangle\\
2\langle \pa_{x_1}\bd^\ve, \pa_{x_2} \bd^\ve \rangle & |\pa_{x_2} \bd^\ve|^2-|\pa_{x_1}\bd^\ve|^2
  \end{matrix}
  \right).
  \label{}
\end{equation}
We can now assume that there exists two real number $\alpha, \beta$ such that 
\begin{align*}
  &\left( \nabla\bd^\ve\odot\nabla\bd^\ve-\frac{1}{2}|\nabla\bd^\ve|^2 \mathbb{I}_2 \right)dx\\
  &\rightharpoonup \left( \frac{1}{2}\nabla\d\odot\nabla\d-\frac{1}{2}|\nabla\d|^2 \mathbb{I}_2 \right)dx+\left( 
  \begin{matrix}
    \alpha&\beta\\
    \beta&-\alpha
  \end{matrix}
  \right)\delta_{(0,0)}
\end{align*}
as convergence of Radon measures. \eqref{eqn:convofEricksen} is true if we can show 
\begin{equation*}
  \alpha=\beta=0.
\end{equation*}
We apply the same Pohozaev argument as that in 
\cite{DuHuangWang2020compactness}. Set $\tau^\ve$, $e_\ve$ to be  
\begin{equation}
  \Delta \bd^\ve-\f_\ve(\bd^\ve)=:\tau^\ve
  \label{eqn:bdeq}
\end{equation}
and
\begin{equation*}
  e_\ve(\bd^\ve):=\frac{1}{2}|\nabla\bd^\ve|^2+F_\ve(\bd^\ve). 
\end{equation*}
For any $X\in C^\infty(D, \R^2)$, multiplying \eqref{eqn:bdeq} by $X\cdot \nabla\bd^\ve$ and integrating over $B_r(0)$ we get
\begin{equation}
  \begin{split}
    &\int_{\pa B_r(0)}\langle X\cdot \nabla \bd^\ve, \frac{ x}{|x|} \rangle d\sigma-\int_{B_r(0)}\langle \nabla\bd^\ve\odot\nabla\bd^\ve, \nabla X \rangle dx\\
    &+\int_{B_r(0)}\dv X e_\ve(\bd^\ve) dx
    -\int_{\pa B_r(0)}e_\ve(\bd^\ve)\langle X, \frac{\x}{|\x|} \rangle d\sigma\\
    &=\int_{B_r(0)}\langle X\cdot \nabla\bd^\ve, \tau^\ve \rangle dx.
  \end{split}
  \label{eqn:PohozaevX}
\end{equation}
If we choose $X(\x)=\x$, then \eqref{eqn:PohozaevX} becomes
\begin{align*}
  &r\int_{\pa B_r(0)}\left|\frac{\pa\bd^\ve}{\pa r}\right|^2d\sigma+\int_{B_r(0)} 2 F_\ve(\bd^\ve)dx\\
  &-r\int_{\pa B_r(0)}e_\ve(\bd^\ve)d\sigma=\int_{B_r(0)}|\x|\langle\frac{\pa\bd^\ve}{\pa r}, \tau^\ve\rangle dx. 
\end{align*}
Hence 
\begin{align*}
  &\int_{\pa B_r(0)}e_\ve(\bd^\ve)d\sigma=\int_{\pa B_r(0)}\left|\frac{\pa \bd^\ve}{\pa r}\right|d\sigma\\
  &+\frac{1}{r}\int_{B_r(0)} 2F_\ve(\bd^\ve)dx-\frac{1}{r}\int_{B_r(0)}|\x|\langle \frac{\pa \bd^\ve}{\pa r}, \tau^\ve \rangle dx.
\end{align*}
Integrating from $r$ to $R$ yields 
\begin{equation}
  \begin{split}
   & \int_{B_R(0)\setminus B_r(0)}e_\ve(\bd^\ve)dx=\int_{B_R(0)\setminus B_r(0)}\left|\frac{\pa\bd^\ve}{\pa r}\right|^2dx\\
    &+\int_{r}^{R}\frac{1}{\tau}\int_{B_\tau(0)}\left( 2F_\ve(\bd^\ve)-|\x|\langle \frac{\pa\bd^\ve}{\pa r}, \tau^\ve \rangle  \right) dxd\tau.
  \end{split}\label{eqn:rR}
\end{equation}
Since $\Sigma_t={(0,0)}$, then there exists $\gamma>0$ such that 
\begin{equation*}
  e_\ve(\bd^\ve)d\x\rightharpoonup \frac{1}{2}|\nabla \d|^2 d\x+\gamma \delta_{(0,0)} 
\end{equation*}
as convergence of Radon measure. By sending $\ve\to 0$ in \eqref{eqn:rR} we get
\begin{equation}
	\begin{split}
  &\int_{B_R(0)\setminus B_r(0)}\frac{1}{2}|\nabla\d|^2 dx\\
  &\ge \int_{B_R(0)\setminus B_r(0)}\left|\frac{\pa \d}{\pa r}\right|^2 dx+\int_{r}^{R}\frac{1}{\tau}\liminf_{\ve\to 0}\int_{B_\tau(0)} 2F_\ve(\bd^\ve)dxd\tau\\
  &+\liminf_{\ve\to 0}\int_r^R -\frac{1}{\tau}\int_{B_\tau(0)}|\x|\langle \frac{\pa \bd^\ve}{\pa r}, \tau^\ve \rangle dxd\tau.
\end{split}
  \label{eqn:Rr2}
\end{equation}
Notice that
\begin{align*}
  &\left|\int_{r}^{R}-\frac{1}{\tau}\liminf_{\ve\to0}\int_{B_\tau(0)}|\x|\langle \frac{\pa\bd^\ve}{\pa r}, \tau^\ve \rangle dx d\tau\right|\\
  &\le \limsup_{\ve\to0} \int_{0}^{R} \left\|\tau^\ve\right\|_{L^2(B_\tau(0))}\left\|\nabla\bd^\ve\right\|_{L^2(B_\tau(0))}d\tau\\
  &=\mathcal{O}(R). 
\end{align*}
As a consequence, we claim that 
\begin{equation}
  2F_\ve(\bd^\ve)\to 0 \text{ in } L^1(B_\delta).\label{eqn:Fconv} 
\end{equation}
For, otherwise, then there exists $\kappa>0$ such that 
\begin{equation*}
  2F_\ve(\bd^\ve)d\x\rightharpoonup \kappa \delta_{(0,0)}.
\end{equation*}
This implies
\begin{equation*}
  \lim_{r\downarrow 0}\int_{r}^{R}\frac{1}{\tau}\liminf_{\ve\to 0}\int_{B_\tau(0)}2F_\ve(\bd^\ve)dxd\tau =\lim_{r\downarrow0} \int_{r}^{R} \frac{\kappa}{\tau}d\tau=\infty.
\end{equation*}
If we choose $X(\x)=(x_1, 0)$ in \eqref{eqn:PohozaevX}, we obtain that 
\begin{equation}
  \begin{split}
    &\frac{1}{2}\int_{B_r(0)}\left( |\pa_{x_2}\bd^\ve|^2-|\pa_{x_1} \bd^\ve|^2 \right)dx+\int_{B_r(0)}F_\ve(\bd^\ve)dx\\
    &=\int_{B_r(0)}x_1\langle \pa_{x_1}\bd^\ve, \tau^\ve \rangle dx+\int_{\pa B_r(0)}\frac{x_1^2}{r}e_\ve(\bd^\ve)d\sigma\\
    &-\int_{\pa B_r(0)}x_1\langle \pa_{x_1} \bd^\ve, \frac{\pa \bd^\ve}{\pa r} \rangle d\sigma.
  \end{split}
  \label{eqn:Pohox}
\end{equation}
Since  $e_\ve(\bd^\ve)d\x\rightharpoonup \frac{1}{2}|\nabla\d|^2 d\x$ in $B_{2r}\setminus B_{\frac{r}{2}}$ for $r>0$, it is easy to see
\begin{align*}
  &\int_{\pa B_r(0)}x_1\langle \pa_{x_1} \bd^\ve, \frac{\pa \bd^\ve}{\pa r} \rangle d\sigma\to \int_{\pa B_r(0)}x_1\langle \pa_{x_1}\d, \frac{\pa\d}{\pa r} \rangle d\sigma, \\
  &\int_{\pa B_r(0)}\frac{x_1^2}{r}e_\ve(\bd^\ve)d\sigma\to \frac{1}{2}\int_{\pa B_r}\frac{x_1^2}{r}|\nabla\d|^2d\sigma, 
\end{align*}
and by \eqref{eqn:Fconv}, 
\begin{equation*}
  \int_{B_r(0)}F_\ve(\bd^\ve) dx\to 0.
\end{equation*}
With the fact that 
\begin{equation*}
  \left|\int_{B_r}x_1\langle \pa_{x_1}\bd^\ve, \tau^\ve \rangle dx\right|=\mathcal{O}(r),
\end{equation*}
by sending $\ve\to 0$ in \eqref{eqn:Pohox} we obtain
$$
\frac{1}{2}\int_{B_r(0)}\left( |\pa_{x_2}\d|^2-|\pa_{\x_1}\d|^2 \right)dx+\alpha=\mathcal{O}(r)$$
which implies $\alpha=0$ after sending $r\to 0$.

Similarly, if we choose $X(\x)=(0,x_1)$ in \eqref{eqn:PohozaevX}, by performing the same argument we will arrive at
$$  \frac{1}{2}\int_{B_r(0)}\langle \pa_{x_1}\d, \pa_{x_2} \d \rangle dx +\beta=\mathcal{O}(r).
 $$
 Hence $\beta=0$. This implies almost surely convergence of Ericksen stress tensor field \eqref{eqn:convofEricksen}.  From \eqref{eqn:unifpnorm} and \eqref{eqn:law} we can conclude that for any $1<p<\infty$, it holds
 \begin{align}\label{eqn:uniforminteg}
	&\E' \sup_{0\le t\le T}\left[ \|\bu^\ve(t)\|^2+\|\nabla\bd^\ve(t)\|^2+\int_{D\times\{t\}} F_\ve(\bd^\ve)dx \right]^p \\
	& +\E'\left[\int_{0}^{T}(\|\nabla\bu^\ve\|^2+\|\Delta\bd^\ve-\f_\ve(\bd^\ve)\|^2)dt\right]^p\nonumber\\
	&=\E \sup_{0\le t\le T}\left[ \|\u^\ve(t)\|^2+\|\nabla\d^\ve(t)\|^2+\int_{D\times\{t\}} F_\ve(\d^\ve) dx\right]^p \nonumber\\
	& +\E\left[\int_{0}^{T}(\|\nabla\u^\ve\|^2+\|\Delta\d^\ve-\f_\ve(\d^\ve)\|^2)dt\right]^p\nonumber\\
	&\le C. \nonumber
\end{align}
  Thus we have for any $\xi\in L^2(\Omega'; {\bf J})$, it holds
 \begin{align}\label{eqn:converofu}
   &\quad\lim_{\ve\to0}\E'\left[\int_D \langle M_{\bu^\ve}(t), \xi \rangle dx\right]\\
   &=\lim_{\ve\to 0} \E'\Big[ \int_D \langle \bu^\ve(t)-\u_0^\ve, \xi \rangle dx\nonumber\\
   &\quad+\int_{0}^{t}\int_D(\langle \nabla \bu^\ve, \nabla \xi \rangle-\langle  \bu^\ve\otimes\bu^\ve+\nabla  \bd^\ve\odot \nabla\bd^\ve-\frac{1}{2}|\nabla \bd^\ve|^2\mathbb{I}_2, \nabla \xi \rangle) dxds 
   \Big]\nonumber\\
   &=\E'\Big[ \int_D \langle \u(t)-\u_0, \xi \rangle dx\nonumber\\
   &\quad+\int_{0}^{t}\int_D(\langle \nabla \u, \nabla \xi \rangle-\langle \u\otimes\u+\nabla\d\odot\nabla\d-\frac{1}{2}|\nabla \d|^2\mathbb{I}_2, \nabla \xi \rangle) dxds \Big]\nonumber\\
   &=\E'\left[ \int_D\langle M_{\u}(t), \xi \rangle dx \right]. \nonumber
 \end{align}

Now we turn to the convergence of $\bd^\ve$. 
We claim that up to a subsequence, 
\begin{equation}
  \Delta \bd^\ve-\f_\ve(\bd^\ve)\rightharpoonup \Delta \d+|\nabla\d|^2\d \text{ in }L^2(\Omega'\times[0,T]\times D).
  \label{eqn:convfd}
\end{equation}
From \eqref{eqn:unifEnergy} we can assume that there exists ${\bf g}\in L^2(\Omega'\times[0,T]\times D)$ such that 
\begin{equation*}
  \Delta\bd^\ve-\f_\ve(\bd^\ve)\rightharpoonup {\bf g} \text{ in }L^2(\Omega'\times[0,T]\times D).
\end{equation*}
First we claim that
\begin{equation}
  {\bf g}\perp \d \text{ for almost all }(\omega', t, x)\in \Omega'\times[0,T]\times D
  \label{eqn:gperp}.
\end{equation}
In fact, for any test function $\phi=\phi(\omega',x)$, if we apply the It\^{o} formula to 
$$\Psi(\bd^
\ve)=\int_D \frac{|\bd^\ve|^2}{2}\phi dx,$$
 it hold that (see Appendix A) 
\begin{align*}
  &\E'\left[ \int_D \frac{|\bd^\ve|^2(t)}{2}\phi dx \right]-\E'\left[ \int_D \frac{|\bd^\ve|^2(t-\delta)}{2}\phi dx\right]\\
  &=-\E'\left[\int_{t-\delta}^{t}\int_D\phi \bu^\ve\cdot \nabla \frac{|\bd^\ve|^2}{2}dxds\right]+\E'\left[\int_{t-\delta}^{t}\int_D\langle \Delta\bd^\ve-\f_\ve(\bd^\ve), \bd^\ve \rangle \phi dxds\right].
\end{align*}
Now we pass $\ve$ to 0, using the fact that $|\d|=1$ for almost all $(\omega',t,x)\in \Omega'\times[0,T]\times D$ we get
\begin{equation}
  \E'\left[ \int_{t-\delta}^{t}\int_D \langle {\bf g}, \d \rangle \phi dxds\right]=0.
  \label{}
\end{equation}
Since $\phi$ and $\delta$ can be arbitrary,  $\langle {\bf g}, \d \rangle=0$ for almost all $(\omega',t,x)\in \Omega'\times[0,T]\times D$. Hence \eqref{eqn:gperp} holds. 
By taking the cross product of \eqref{eqn:convfd} with $\bd^\ve \phi$ we get
\begin{align*}
  0&=\lim_{\ve\to 0} \E'\left[ \int_{0}^{T}\int_D \langle (\Delta \bd^\ve -\f_\ve(\bd^\ve)-{\bf g})\times \bd^\ve,  \phi\rangle dxdt \right]\\
  &=\lim_{\ve\to0}\E' \left[ \int_{0}^{T}\int_D \langle \nabla\cdot (\nabla \bd^\ve\times\bd^\ve), \phi \rangle dxdt-\int_{0}^{T}\int_D \langle {\bf g}\times\bd^\ve, \phi \rangle dxdt\right]\\
  &=\lim_{\ve\to 0}\E'\left[ -\int_{0}^{T}\int_D \langle \nabla\bd^\ve\times\bd^\ve, \nabla \phi \rangle dxdt-\int_{0}^{T}\int_D \langle {\bf g}\times \bd^\ve, \phi\rangle dxdt\right]
  \\
  &=\E'\left[ -\int_{0}^{T}\int_D \langle \nabla \d\times\d, \nabla \phi \rangle dxdt-\int_{0}^{T}\int_D \langle {\bf g}\times \d, \phi \rangle dxdt \right].
\end{align*}
This implies $({\bf g}-\Delta \d)\times \d=0$ and hence there exists $\lambda=\lambda(\omega',t, x):\Omega'\times[0,T]\times D\to \R$ such that 
\begin{equation*}
  {\bf g}-\Delta \d=\lambda \d.
\end{equation*}
From \eqref{eqn:gperp} and $\langle \Delta \d, \d\rangle=-|\nabla \d|^2 \d$ we 
get
\begin{equation*}
  \lambda=\langle {\bf g}-\Delta \d, \d \rangle=|\nabla \d|^2.
\end{equation*}
Thus \eqref{eqn:convfd} holds. Thanks to \eqref{eqn:conv}, \eqref{eqn:uniforminteg} and \eqref{eqn:convfd}, we have for any $\zeta\in L^2(\Omega'; H^1(D,\R^3))$, it holds
\begin{align}\label{eqn:converofd}
  &\quad \lim_{\ve\to0}\E'\left[ \int_D\langle M_{\bd^\ve}(t), \zeta \rangle dx\right]\\
  &=\lim_{\ve\to0}\E'\Big[ \int_{D}\langle \bd^\ve(t)-\d_0^\ve, \zeta \rangle dx\nonumber \\
  &-\int_{0}^{t}\int_D(\langle \bu^\ve\otimes\bd^\ve, \nabla \zeta \rangle-\langle \Delta\bd^\ve-\f_\ve(\bd^\ve), \zeta \rangle) dxds\Big]\nonumber\\
  &-\lim_{\ve\to0}\E'\left[\frac{1}{2}\int_{0}^{t}\int_D \langle (\bd^\ve\times\h)\times \h, \zeta \rangle dxds\right]\nonumber\\
  &=\E'\left[ \int_D \langle \d(t)-\d_0, \zeta \rangle dx-\int_{0}^{t}\int_D (\langle \u\otimes\d, \nabla \zeta \rangle- \langle \Delta\d+|\nabla \d|^2\d, \zeta \rangle) dxds \right]\nonumber\\
  &\qquad-\E'\left[ \frac{1}{2}\int_{0}^{t}\int_D \langle (\d\times\h)\times\h, \zeta \rangle dxds\right]\nonumber\\
  &=\E'\left[ \int_D \langle M_{\d}(t), \zeta \rangle dx\right]. \nonumber
\end{align}
Taking the limit $\ve\to 0$ in \eqref{eqn:globeng} and applying the lower semicontinuity yields
\eqref{eqn:defenergy}. 
 
To finish the construction, we need to show that for every $t\in(0,T]$
\begin{align}
  \int_{0}^{t} S(\bu^\ve)d\bW_1^\ve(s)&\to \int_{0}^{t} S(\u) dW'_1(s) \text{ in }L^2(\Omega; L^2(D)), \label{3.23}\\
   \int_{0}^{t}(\bd^\ve \times\h)d\bW_2^\ve(s)&\to \int_{0}^{t}(\d\times\h)dW_2'(s) \text{ in } L^2(\Omega;L^2(D)). 
  \label{3.24}
\end{align}
For this purpose, we adapt the strategy from \cite{brzezniak2019some}. Let $\mathcal{N}$ be the set of null sets of $\F'$ and for any $t\ge 0$ and $\ve>0$, let
\begin{align*}
  \hat{\F}_t^\ve&:=\sigma\left( \sigma\left( (\bu^\ve(s), \bd^\ve(s), \bW_1^\ve(s), \bW_2^\ve(s)); s\le t \right)\cup \mathcal{N} \right), \\
  \F'_t&:=\sigma\left( \sigma\left( (\u(s), \d(s), W'_1(s), W'_2(s)); s\le t \right)\cup \mathcal{N} \right). 
\end{align*}
Since $\mathcal{L}(\bu^\ve, \bd^\ve, \bW_1^\ve, \bW_2^\ve)=\mathcal{L}(\u^\ve, \d^\ve, W_1, W_2)$, $(\bW_1^\ve, \bW_2^\ve)$ form a sequence of  cylindrical Wiener processes. Moreover, for $0\le s<t\le T$ the increments $(\bW_1^\ve(t)-\bW_1^\ve(s), \bW_2^\ve(t)-\bW_2^\ve(s))$ are independent of $\hat{\mathcal{F}}_r^\ve$ for $r\in[0,s]$. Let $k\in \mathbb{N}$ and $s_0=0<s_1<\cdots<s_k\le T$ be a partition of $[0,T]$. By the characterization of $K_2$-valued $K_1$-cylindrical Wiener process \cite[Remark 2.8]{brzezniak2019some},  for each $\xi\in K_2^*$ we have
\begin{eqnarray*}
  \E'\left[ e^{i \sum_{j=1}^{k}\left\langle \xi, \bW_1^\ve(s_j)-\bW_1^\ve(s_{j-1}) \right\rangle_{K_2^*,K_2}} \right]&=& \E\left[ e^{i \sum_{j=1}^{k}\left\langle \xi, W_1(s_j)-W_1(s_{j-1}) \right\rangle_{K_2^*, K_2}} \right]\\
  &=&e^{-\frac{1}{2}\sum_{j=1}^{k}(s_j-s_{j-1})|\xi|_{K_1}^2}. 
\end{eqnarray*}
Thanks to \eqref{eqn:conv} and the Lebesgue Dominated Convergence Theorem, we have
\begin{align*}
  \lim_{\ve\to0}\E'\left[ e^{i\sum_{j=1}^{k}\left\langle \xi, \bW_1^\ve(s_j)-\bW_1^\ve(s_{j-1}) \right\rangle_{K_2^*, K_2}} \right]&=\E'\left[ e^{i\sum_{j=1}^{k}\left\langle \xi, W_1'(s_1)-W_1'(s_{j-1}) \right\rangle_{K_2^*, K_2}} \right]\\
  &=e^{-\frac{1}{2}\sum_{j=1}^{k}(s_j-s_{j-1})|\xi|_{K_1}^2}. 
\end{align*}
Hence the finite dimensional distribution of $W_1'$ is Gaussian. The same argument also works for $W_2'$. Next we want to show that $(W_1'(t)-W_1'(s), W_2'(t)-W_2'(s)), 0\le s<t\le T$ is independent of $\mathcal{F}'_r$ for $r\in [0,s]$. Consider $\{\phi_j\}_{j=1}^k\in C_b(W^{-2,\tilde{p}}(D)\times L^{\frac{4}{3}}(D)), \{\psi_j\}_{j=1}^k\in C_b(K_2\times\R)$, let $0\le t_1<\cdots<r_k\le s<t\le T$, $\psi\in C_b(K_2), \zeta\in C_b(\R)$. 
\begin{align}\label{eqn:indep}
  &\E'\Bigg[ \left(\prod_{j=1}^k \phi_j(\bu^\ve(r_j), \bd^\ve(r_j))\prod_{j=1}^k\psi_j(\bW^\ve_1(r_j), \bW^\ve_2(r_j))\right)\\
  &\qquad \times \psi(\bW_1^\ve(t)-\bW_1^\ve(s))\zeta(\bW_2^\ve(t)-\bW_2^\ve(s))\Bigg]\nonumber\\
  =&\E'\left[ \prod_{j=1}^k \phi_j(\bd^\ve(r_j), \bd^\ve(r_j))\prod_{j=1}^k \psi_j(\bW_1^\ve(r_j), \bW_2^\ve(r_j)) \right]\nonumber\\
  &\qquad \times \E'\left[ \psi(\bW_1^\ve(t)-\bW_1^\ve(s)) \right]\E'\left[ \zeta(\bW_2^\ve(t)-\bW_2^\ve(s)) \right].\nonumber
\end{align}
Again by the Lebesgue Dominated Convergence theorem, if we send $\ve\to 0$ in \eqref{eqn:indep} we can see \eqref{eqn:indep}  also holds for $(\u, \d, W_1', W_2')$ in the limit. Furthermore, it is easy to show that $W_1'$ is independent of $W_2'$.  

For any $\delta>0$, let $\eta_\delta$ be a standard mollifier with support in $(0,t)$. Define
\begin{equation*}
  S^\delta(\u(s))=\int_{-\infty}^{\infty}\eta_\delta(s-r)S(\u(r)dr. 
\end{equation*}
Let $M_{\bu^\ve}^\delta$ and $M_{\u}^\delta$ be respectively defined by 
\begin{align*}
  M_{\bu^\ve}^\delta(t)&=\int_{0}^{t}S^\delta(\bu^\ve(s))d\bW_1^\ve(s), \\
  M_{\u}^\delta(t)&=\int_{0}^{t}S^\delta(\u(s))dW_1'(s).
\end{align*}
By the property of mollifiers, we can get for any ${\bf v}\in {\bf H}$
\begin{equation*}
  \lim_{\delta\to0}\E'\int_{0}^{t}\|S^\delta({\bf v}(s))-S({\bf v}(s))\|_{\mathcal{L}_2(K_1, {\bf H})}^2ds=0.
\end{equation*}
Hence, for any $t\in(0,T]$, we have the following uniform approximation
\begin{equation}\label{eqn:deltaapproximate}
  \lim_{\delta\to 0}\sup_{0<\ve<1}\E'\left\|M^\delta_{\bu^\ve}(t)-\int_0^t S(\bu^\ve)d\bW_1^\ve(s)\right\|^2=0,
\end{equation}
and 
\begin{equation}\label{eqn:delta2}
	\lim_{\delta\to 0}\E'\left\|M_{\u}^\delta(t)-\int_0^t S(\u)dW_1'(s)\right\|^2=0. 
\end{equation}
Next, we need to show that for any $\delta>0$
\begin{equation}
  \lim_{\ve\to0}\E'\left\|M^\delta_{\bu^\ve}(t)-M^\delta_{\u}(t)\right\|^2=0.
  \label{eqn:deltaepsilonconv}
\end{equation}
If we write $\bW_1^\ve(t)=\sum_{i=1}^{\infty} \bB_i^\ve(t)e_i$ and $W_1'(t)=\sum_{i=1}^{\infty}B_i'(t)e_i$, where $\{\bB_i^\ve\}_{i=1}^\infty$, $\{B_i'(t)\}_{i=1}^\infty$ are i.i.d. stardard Brownian motions, then 
\begin{equation*}
  M_{\bu^\ve}^\delta(t)-M_\u^\delta(t)=\sum_{i=1}^{\infty}\int_{0}^{t}S^\delta(\bu^\ve(s))(e_i)d\bB_i^\ve(s)-\sum_{i=1}^{\infty}\int_{0}^{t}S^\delta(\u(s))(e_i)dB_i'(s).
\end{equation*}
By Young's convolution inequality, we have that
\begin{equation*}
  \E'\int_{0}^{t}\| S^\delta(\u(s))\|_{\mathcal{L}_2(K_1, {\bf H})}^2 ds\le  C \E'\int_{0}^{t}\|S(\u(s))\|_{\mathcal{L}_2(K_1, {\bf H)}}^2 ds\le C. 
\end{equation*}
Thus, for any $\gamma>0$, there exists an $N\in \mathbb{N}_+$ such that 
\begin{equation*}
  \sum_{i=N+1}^{\infty} \E'\int_{0}^{t}\| S^\delta(\u(s))(e_i)\|^2 ds<\gamma.
\end{equation*}
Since
\begin{equation*}
  \lim_{\ve\to0}\E'\int_{0}^{t}\|S^\delta(\bu^\ve(s))-S^\delta(\u(s))\|_{\mathcal{L}_2(K_1, {\bf H})}^2 ds=0,
\end{equation*}
there exists an $\ve_0>0$ such that for $0<\ve<\ve_0$, 
\begin{equation*}
  \sum_{i=N+1}^{\infty}\E'\int_{0}^{t}\|S^\delta(\bu^\ve(s))(e_i)\|^2 ds<2\gamma.
\end{equation*}
Now we split $M_{\bu^\ve}^\delta(t)-M_{\u}^\delta(t)$ into three parts
\begin{align*}
  &M^\delta_{\bu^\ve}(t)-M_{\u}^\delta(t)=\sum_{i=1}^{N}\left( \int_{0}^{t}S^\delta (\bu^\ve(s))(e_i)d\bB_i^\ve(s)-\int_{0}^{t}S^\delta(\u(s))(e_i)dB_i'(s) \right)\\
  &+\sum_{i=N+1}^{\infty}\int_{0}^{t}S^\delta(\bu^\ve(s))(e_i)d\bB_i^\ve(s)\\
  &+\sum_{i=N+1}^{\infty}\int_{0}^{t}S^\delta(\u(s))(e_i)dB_i'(s):=J^\delta_{\ve, 1}(t)+J^\delta_{\ve, 2}(t)+J^{\delta}_{\ve, 3}(t). 
\end{align*}
By the It\^{o} isometry, we have that
\begin{align*}
  &\E'\|J_{\ve, 2}^\delta\|^2=\sum_{i=N+1}^{\infty}\E'\int_{0}^{t}\|S^\delta(\bu^\ve(s))(e_i)\|^2 ds<2\gamma,\\
  &\E'\|J_{\ve, 3}^\delta\|^2=\sum_{i=N+1}^{\infty}\E'\int_{0}^{t}\|S^\delta(\u(s))(e_i)\|^2 ds<\gamma.
\end{align*}
For $J_{\ve, 1}(t)$, we write
\begin{align*}
  &J_{\ve, 1}^\delta(t)=\sum_{i=1}^{N}\left( \int_{0}^{t}S^\delta(\bu^\ve(s))(e_i)d\bB_i^\ve(s)-\int_{0}^{t}S^\delta(\bu^\ve(s))(e_i)dB_i'(s) \right)\\
  &+\sum_{i=1}^{N}\left( \int_{0}^{t}S^\delta(\bu^\ve(s))(e_i)dB'_i(s)-\int_{0}^{t}S^\delta(\u(s))(e_i)dB_i'(s) \right)\\
  &:=I_{\ve, 1}^\delta+I_{\ve, 2}^\delta.
\end{align*}
For $I^\delta_{\ve, 1}(t)$, by integration by parts we obtain that
\begin{align*}
  I^\delta_{\ve, 1}(t)&=\sum_{i=1}^{N}\left(\int_{0}^{t}[\eta'_\delta\star S(\bu^\ve(s))](e_i)B_i'(s)ds-\int_{0}^{t}[\eta'_\delta\star S(\bu^\ve(s))](e_i)\bB^\ve_i(s)ds\right)\\
  &=-\sum_{i=1}^{N}\left(\int_{0}^{t}[\eta'_\delta\star S(\bu^\ve(s))](e_i)[\bB_i^\ve(s)-B_i'(s)]ds\right).
\end{align*}
From the Burkholder--Davis--Gundy inequality, we get for any $p>1$, any $i=1,2,...,N$,
\begin{equation}
  \sup_{\ve>0}\E'\sup_{s\in[0,T]}\left( |\bB^\ve_i(s)|^p+|B_i'(s)|^p \right)\le CT^{\frac{p}{2}}.\label{eqn:uninteB}
\end{equation}
Hence, by the uniform integrability \eqref{eqn:uninteB} and the almost surely convergence \eqref{eqn:conv} 
we have that for $i=1, 2, \dots, N$, 
\begin{equation*}
  \lim_{\ve\to0}\E'\int_{0}^{t}|\bB_i^\ve(s)-B'_i(s)|^p ds=0.
\end{equation*}
This implies 
\begin{align*}
  &\E'\|I_{\ve, 1}^\delta(t)\|^2=\E'\left\|\sum_{i=1}^{N}\int_{0}^{t}[\eta_\delta'\star S(\bu^\ve(s))](e_i)(\bB_i^\ve(s)-B_i'(s))ds\right\|^2\\
  &\le N\sum_{i=1}^{N} \E'\left\|\int_{0}^{t}[\eta_\delta'\star S(\bu^\ve(s))](e_i)(\bB_i^\ve(s)-B_i'(s))ds\right\|^2\\
  &\le N\sum_{i=1}^{N}\E'\left[ \int_{0}^{t}\|[\eta_\delta'\star S(\bu^\ve(s))](e_i)\||\bB_i^\ve(s)-B_i'(s)|ds \right]^2\\
  &\le N\sum_{i=1}^{N}\E'\left[ \int_{0}^{t}\|\eta_\delta'\star S(\bu^\ve(s))(e_i)\|^2 ds \int_{0}^{t}|\bB_i^\ve(s)-B_i'(s)|^2ds \right]\\
  &\le \frac{CN}{\delta^2}\sum_{i=1}^{N} \E'\left[ \int_{0}^{t}\|S(\bu^\ve)(s)\|_{\mathcal{L}_2(K_1, {\bf H})}^2 ds \int_{0}^{t}|\bB_i^\ve(s)-B_1'(s)|^2ds \right]\\
  &\le \frac{CN}{\delta^2}\sum_{i=1}^{N}\E'\left[ \int_{0}^{t}(1+\|\bu^\ve\|^2)ds\int_{0}^{t}|\bB_i^\ve(s)-B_i'(s)|^2 ds \right]\\
  &\le \frac{CNT^{\frac{3}{2}}}{\delta^2}\left( \E'\sup_{0\le s\le t}(1+\|\bu^\ve(s)\|)^4 \right)^{\frac{1}{2}}\left( \E'\int_{0}^{t}|\bB^\ve_i(s)-B_i'(s)|^4 ds \right)^{\frac{1}{2}}\\
  &\le \frac{CNT^{\frac{3}{2}}}{\delta^2}\sum_{i=1}^{N}\left( \E'\int_{0}^{t}|\bB_i^\ve(s)-B'_i(s)|^4 ds \right)^{\frac{1}{2}}\to 0,
\end{align*}
as $\ve\to0$.
Using a similar argument, we can show that 
\begin{equation*}
  \lim_{\ve\to0}\E'\|I_{\ve, 2}^\delta(t)\|^2=0.
\end{equation*}
Since $\gamma$ can be arbitrarily small,  we get
\begin{equation*}
  \lim_{\ve\to0}\E'\left[ \|J_{\ve,1}^\delta\|^2+\|J_{\ve,2}^\delta(t)\|^2+\|J_{\ve,3}^\delta(t)\|^2 \right]=0, \ \forall t\in(0,T], 
\end{equation*}
This implies \eqref{eqn:deltaepsilonconv}.  Then we can conclude from  \eqref{eqn:deltaapproximate}, \eqref{eqn:delta2} and \eqref{eqn:deltaepsilonconv}  that for every $t\in(0,T]$, 
\begin{equation*}
  \lim_{\ve\to0}\E'\left\|\int_0^tS(\bu^\ve(s))d\bW_1^\ve(s)-\int_0^t S(\u(s))dW_1'(s)\right\|^2=0. 
\end{equation*}
Similarly, we can show
\begin{equation*}
  \lim_{\ve\to0}\E'\left\|\int_0^t(\bd^\ve\times\h)d\bW_2^\ve(s)-\int_{0}^{t}(\d\times
  \h)dW_2'(s)\right\|^2=0. 
\end{equation*}
Hence, the convergence of martingale terms \eqref{3.23} and \eqref{3.24} holds. Putting \eqref{eqn:converofu}, \eqref{eqn:converofd}, \eqref{3.23} and \eqref{3.24} together completes the proof.

\appendix
\section{It\^{o}'s formulas for functionals of $\d$}
Consider the functional 
\begin{equation*}
	\Psi(\d^\ve):=\int_D \frac{|\d^\ve|^2}{2}\phi dx.
\end{equation*}
It is easy to obtain the first and and second Fr\'{e}chet derivatives of $\Psi(\d^
\ve)$
\begin{align*}
	\Psi'(\d^\ve)[{\bf g}]&=\int_D \langle \d^\ve, {\bf g} \rangle \phi dx, \\
	\Psi''(\d^\ve)[{\bf g}, {\bf g}]&=\int_D \langle {\bf g}, {\bf g} \rangle \phi dx.
\end{align*}
Applying the It\^{o} formula to $\Psi(\d^\ve)$ gives
\begin{align*}
	d\Psi(\d^\ve)&=\Psi'(\d^\ve)[d\d^\ve]+\frac{1}{2}\Psi_\ve''(\d^\ve)[d\d^\ve, d\d^\ve].
\end{align*}
Since, 
\begin{equation*}
	d \d^\ve=\underbrace{(-\u^\ve\cdot \nabla\d^\ve+\Delta\d^\ve-\f_\ve(\d^\ve)+\frac{1}{2}(\d^\ve\times\h)\times\h)}_{ {\bf j}}dt+\underbrace{(\d^\ve\times\h)}_{ {\bf k}}dW_2, 
\end{equation*}
we then obtain that for $0<\delta<t$, 
\begin{align*}
	&\Psi(\d^\ve)(t)-\Psi(\d^\ve)(t-\delta)=\int_{t-\delta}^{t}\left(\Psi'(\d^\ve)[{\bf j}(s)]+\frac{1}{2}\Psi''(\d^\ve)[{\bf k}(s), {\bf k}(s)]\right)ds\\
	&+\int_{t-\delta}^{t}\Psi'(\d^\ve)[{\bf k}(s)]dW_2(s)\\
	&=\int_{t-\delta}^{t}\int_D \langle -\u^\ve\cdot \nabla\d^\ve, \d^\ve \rangle\phi dxds+\int_{t-\delta}^{t}\int_D \langle \Delta\d^\ve-\f_\ve(\d^\ve), \d^\ve \rangle\phi dxds\\
	&+\frac{1}{2}\int_{t-\delta}^{t}\int_D \langle (\d^\ve\times\h)\times\h, \d^\ve \rangle\phi dxds+\frac{1}{2}\int_{t-\delta}^{t}\int_D |\d^\ve\times\h|^2 \phi dxds\\
	& +\int_{t-\delta}^{t}\int_D \langle \d^\ve\times\h, \d^\ve \rangle\phi dxdW_2(s)\\
	&=\int_{t-\delta}^{t}-\u^\ve\cdot\nabla  \frac{|\d^\ve|^2}{2}\phi dxds+\int_{t-\delta}^{t}\int_D \langle \Delta\d^\ve-\f_\ve(\d^\ve), \d^\ve \rangle\phi dxds,
\end{align*}
where we use the fact  the vector triple product
\begin{equation}
	\langle (\d^\ve\times\h)\times\h, \d^\ve \rangle=-|\d^\ve\times\h|^2.\label{eqn:tripleprod}
\end{equation}
and
\begin{equation}
	\langle \d^\ve\times\h, \d^\ve \rangle =0.\label{eqn:perp}
\end{equation}

Recall the energy functional 
\begin{equation*}
  \Phi_\ve(\d^\ve)=\frac{1}{2}\|\nabla \d^\ve\|^2+\int_D F_\ve(\d^\ve)dx.
\end{equation*}
The first and second Fr\'{e}chet derivatives of $\Phi_\ve$ are given by 
\begin{align*}
  \Phi_\ve'(\d^\ve)[{\bf g}] &=\int_D \left(\langle \nabla \d^\ve, \nabla {\bf g} \rangle+ \langle \f_\ve(\d^\ve), {\bf g} \rangle \right)dx\\
  &=\int_D \langle -\Delta\d^\ve+\frac{|\d^\ve|^2-1}{\ve^2}\d^\ve, {\bf g} \rangle dx,
  \\
  \Phi_\ve''(\d^\ve)[{\bf g}, {\bf g}]&=\int_D \left(\langle \nabla {\bf g}, \nabla{\bf g} \rangle+\frac{|\d^\ve|^2-1}{\ve^2}|{\bf g}|^2+\frac{2}{\ve^2}\langle \d^\ve, {\bf g} \rangle^2 \right)dx 
\end{align*}
for every ${\bf g}\in H^1(D;\R^3)$.
Then, the It\^{o} formula for $\Phi_\ve(\d^\ve)$ reads
\begin{align*}
  d\Phi_\ve(\d^\ve)&=\Phi'_\ve(\d^\ve)[d\d^\ve]+\frac{1}{2}\Phi''_\ve(\d^\ve)[d\d^\ve, d\d^\ve].
\end{align*}
From the identity \eqref{eqn:tripleprod} and \eqref{eqn:perp} we obtain
\begin{align*}
  &\Phi_\ve(\d^\ve)(t)-\Phi_\ve(\d_0)\\
  &=\int_{0}^{t}\left(\Phi_\ve'(\d^\ve)[{\bf j}(s)]+\frac{1}{2}\Phi_\ve''(\d^\ve)[{\bf k}(s), {\bf k}(s)]\right)ds+\int_{0}^{t}\Phi_\ve'(\d^\ve
  )[{\bf k}(s)]dW_2(s)\\
  &=\int_{0}^{t}\int_D \langle\u^\ve\cdot \nabla\d^\ve, \Delta\d^\ve-\f_\ve(\d^\ve)\rangle dxds-\int_{0}^{t}\int_D|\Delta\d^\ve-\f_\ve(\d^\ve)|^2 dxds\\&+\frac{1}{2}\int_{0}^{t}\int_D \langle -\Delta\d^\ve+\frac{|\d^\ve|^2-1}{\ve^2}\d^\ve, (\d^\ve\times\h)\times \h \rangle dxds\\
  &+\frac{1}{2}\int_{0}^{t}\int_D \left(|\nabla(\d^\ve\times\h)|^2+\frac{|\d^\ve|^2-1}{\ve^2}|\d^\ve\times\h|^2\right) dxds\\
  &+\frac{1}{2}\int_{0}^{t}\int_D \langle -\Delta\d^\ve+\f_\ve(\d^\ve), \d^\ve\times\h \rangle dxdW_2(s)\\
  &=\int_{0}^{t}\int_D \langle \u^\ve\cdot \nabla\d^\ve, \Delta\d^\ve-\f_\ve(\d^\ve) \rangle dxds-\int_{0}^{t}\int_D |\Delta\d^\ve-\f_\ve(\d^\ve)|^2 dxds\\
  &+\frac{1}{2}\int_{0}^{t}\int_D (\langle \nabla\d^\ve, \nabla( (\d^\ve\times\h)\times\h) \rangle+|\nabla(\d^\ve\times\h)|^2) dxds\\
  &+\frac{1}{2}\int_{0}^{t}\int_D \langle -\Delta\d^\ve+\f_\ve(\d^\ve), \d^\ve\times\h \rangle dx dW_2(s). 
\end{align*}


\begin{thebibliography}{10}
	\bibitem{bensoussan1995stochastic}
	Alain Bensoussan, \emph{Stochastic Navier--Stokes equations}, Acta Applicandae
	Mathematica \textbf{38} (1995), no.~3, 267--304.
	
	\bibitem{bouard2019existence}
	Anne~De Bouard, Antoine Hocquet, and Andreas Prohl, \emph{Existence, uniqueness
		and regularity for the stochastic {E}ricksen--{L}eslie equation}, arXiv:
	1902.05921 (2019).
	
	
	\bibitem{brzezniak2020GL}
	Zdzislaw Brze\'{z}niak, Gabriel Deugou\'{e}, and Paul~Andr\'{e} Razafimandimby, \emph{On
		strong solution to the 2D stochastic Ericksen--Leslie system: A
		Ginzburg--Landau approximation approach}, arXiv preprint arXiv:2011.00100
	(2020).
	
	\bibitem{brzezniak20202d}
	Zdzislaw Brze\'{z}niak, Gabriel Deugou\'{e}, and Paul~Andr\'{e} Razafimandimby, \emph{On the 2D Ericksen-Leslie equations with anisotropic energy and
		external forces}, arXiv preprint arXiv:2005.07659 (2020).
	
	\bibitem{brzezniak2020strong}
	Zdzislaw Brze\'{z}niak, Erika Hausenblas, and Paul ~Andr\'{e} Razafimandimby, \emph{Strong
		solution to stochastic penalised nematic liquid crystals model driven by
		multiplicative Gaussian noise},  arXiv: 2004.00590 (2020).
	
	\bibitem{Brzezniakhausenblas2019NoteonSEL}
	Zdzis{\l}aw Brze\'{z}niak, Erika Hausenblas, and Paul~Andr\'{e} Razafimandimby,
	\emph{A note on the stochastic {E}ricksen-{L}eslie equations for nematic
		liquid crystals}, Discrete Contin. Dyn. Syst. Ser. B \textbf{24} (2019),
	no.~11, 5785--5802.
	
	\bibitem{brzezniak2019some}
	Zdzis{\l}aw Brze{\'z}niak, Erika Hausenblas, and Paul~Andr{\'e} Razafimandimby,
	\emph{Some results on the penalised nematic liquid crystals driven by
		multiplicative noise: weak solution and maximum principle}, Stochastics and
	Partial Differential Equations: Analysis and Computations \textbf{7} (2019),
	no.~3, 417--475.
	
	\bibitem{Brzeniak2019SELwithjumpnoise}
	Zdzis{\l}aw Brze\'{z}niak, Utpal Manna, and Akash~Ashirbad Panda,
	\emph{Martingale solutions of nematic liquid crystals driven by pure jump
		noise in the {M}arcus canonical form}, J. Differential Equations \textbf{266}
	(2019), no.~10, 6204--6283.
	
		
	\bibitem{DipernaMajda1988concentration}
	Ronald~J. DiPerna and Andrew Majda, \emph{Reduced {H}ausdorff dimension and
		concentration-cancellation for two-dimensional incompressible flow}, J. Amer.
	Math. Soc. \textbf{1} (1988), no.~1, 59--95.
	
	\bibitem{DuHuangWang2020compactness}
	Hengrong Du, Tao Huang, and Changyou Wang, \emph{Weak compactness of simplified
		nematic liquid flows in 2D}, arXiv:2006.04210 (2020).
	
	
	\bibitem{ericksen1961conservation}
	Jerald~L. Ericksen, \emph{Conservation laws for liquid crystals}, Transactions
	of the Society of Rheology \textbf{5} (1961), no.~1, 23--34.
	
	\bibitem{Ericksen1962}
	Jerald~L. Ericksen, \emph{Hydrostatic theory of liquid crystals}, Arch. Rational Mech.
	Anal. \textbf{9} (1962), 371--378. \MR{137403}
	
	\bibitem{flandoli1995martsol}
	Franco Flandoli and Dariusz Gatarek, \emph{Martingale and stationary solutions
		for stochastic {N}avier--{S}tokes equations}, Probab. Theory Related Fields
	\textbf{102} (1995), no.~3, 367--391.
	

	
	\bibitem{gyongy1996existence}
	Istv{\'a}n Gy{\"o}ngy and Nicolai Krylov, \emph{Existence of strong solutions
		for it{\^o}'s stochastic equations via approximations}, Probability theory
	and related fields \textbf{105} (1996), no.~2, 143--158.
	
	\bibitem{HuangLinWang2014EricksenLeslieTwoDimension}
	Jinrui Huang, Fanghua Lin, and Changyou Wang, \emph{Regularity and existence of
		global solutions to the {E}ricksen--{L}eslie system in {$\Bbb{R}^2$}}, Comm.
	Math. Phys. \textbf{331} (2014), no.~2, 805--850.
	
	\bibitem{Hong2011}
	Min-Chun Hong, \emph{Global existence of solutions of the simplified Ericksen--Leslie system in dimension two}. 
	Calc. Var. Partial Differential Equations \textbf{40} (2011), no. ~1-2, 15--36. 
	
	 \bibitem{HongXin2011}
	 Min-Chun Hong, and Zhouping Xin,  \emph{Global existence of solutions of the liquid crystal flow for the Oseen--Frank model in $\mathbb R^2$}. 
	 Adv. Math. \textbf{231} (2012), no. ~3-4, 1364--1400.
	
	\bibitem{kortum2019concentrationcancellation}
	Joshua Kortum, \emph{Concentration-cancellation in the ericksen--leslie model},
	Calculus of Variations and Partial Differential Equations \textbf{59} (2020),
	no.~6, 1--16.
	
	\bibitem{Leslie1968}
	Frank~M. Leslie, \emph{Some constitutive equations for liquid crystals}, Arch.
	Rational Mech. Anal. \textbf{28} (1968), no.~4, 265--283.
	
	\bibitem{leslie1992continuum}
	Frank~M. Leslie, \emph{Continuum theory for nematic liquid crystals}, Continuum
	Mechanics and Thermodynamics \textbf{4} (1992), no.~3, 167--175.
	
	\bibitem{lin2010liquid}
	Fang-Hua Lin, Junyu Lin, and Changyou Wang, \emph{Liquid crystal flows in two
		dimensions}, Archive for Rational Mechanics and Analysis \textbf{197} (2010),
	no.~1, 297--336.
	
	\bibitem{lin2010}
	Fang-Hua Lin, and Changyou Wang, \emph{On the uniqueness of heat flow of harmonic maps and hydrodynamic flow of nematic liquid crystals.}
	Chin. Ann. Math. Ser. B \textbf{31} (2010), no. ~ 6, 921--938.
	
	\bibitem{linliu1995nonparabolic}
	Fang-Hua Lin and Chun Liu, \emph{Nonparabolic dissipative systems modeling the
		flow of liquid crystals}, Comm. Pure Appl. Math. \textbf{48} (1995), no.~5,
	501--537.
	
	\bibitem{lin2014recent}
	Fang-Hua Lin and Changyou Wang, \emph{Recent developments of analysis for
		hydrodynamic flow of nematic liquid crystals}, Philosophical Transactions of
	the Royal Society A: Mathematical, Physical and Engineering Sciences
	\textbf{372} (2014), no.~2029, 20130361.
	
	\bibitem{linwang2016weakthreedim}
	Fang-Hua Lin and Changyou Wang, \emph{Global existence of weak solutions of the nematic liquid crystal
		flow in dimension three}, Comm. Pure Appl. Math. \textbf{69} (2016), no.~8,
	1532--1571.
	\bibitem{simon1990sobolev}
	Jacques Simon, \emph{Sobolev, besov and nikolskii fractional spaces: imbeddings
		and comparisons for vector valued spaces on an interval}, Annali di
	Matematica Pura ed Applicata \textbf{157} (1990), no.~1, 117--148.
	\bibitem{Struwe85}
	Michael Struwe, \emph{On the evolution of harmonic mappings of Riemannian surfaces}, Comment. Math. Helvetici, \textbf{60} (1985), 558--581.

	

	
\end{thebibliography}
\end{document}